\documentclass[12pt]{amsart}
\usepackage{fullpage}
\usepackage{colonequals}

\usepackage[english]{babel}   
\usepackage[T1]{fontenc}
\usepackage{lmodern}          

\usepackage[
  protrusion=true,            
  expansion=true,             
  factor=1000,                
  stretch=20,                 
  shrink=20,                  
]{microtype}

\usepackage{amsmath,amssymb,amsthm,verbatim,mathtools,mathrsfs}
\usepackage{xspace}
\usepackage[numbers]{natbib} 
\usepackage[unicode,]{hyperref}
\usepackage{expl3} 
\usepackage{enumitem}
\usepackage{tikz-cd}

\usepackage{xcolor}
\usepackage{comment}

\newtheorem{theorem}{Theorem}[section]
\newtheorem{proposition}[theorem]{Proposition}

\newtheorem{lemma}[theorem]{Lemma}
\newtheorem{corollary}[theorem]{Corollary}

\theoremstyle{definition}
\newtheorem{definition1}[theorem]{Definition}
\newtheorem{remark1}[theorem]{Remark}

\newlist{thmlist}{enumerate}{1}
\setlist[thmlist]{label=(\alph*), ref=\thetheorem(\alph*)}

\usepackage[margin=1.1in]{geometry}
\theoremstyle{remark}

\newcommand{\sgn}{\operatorname{sign}}

\newcommand{\tc}{\operatorname{tc}}
\def\what{\hat}

\newcommand\classic{Permutation Wordle\xspace}
\newcommand\rainbow{Multicolor Permutation Wordle\xspace}
\newcommand\basicstrat{\mbox{\textnormal{\textsc{CircularShift}}}\xspace}
\newcommand{\BC}{\mathbf C}
\newcommand{\Z}{\mathbf Z}
\newcommand{\TP}{\tilde P}
\newcommand{\TZ}{\tilde Z}
\newcommand{\cJ}{\mathcal J}
\newcommand{\K}{\mathcal K}
\newcommand{\Kplus}{{{\mathcal K}^+}}
\newcommand{\FF}{\mathcal F}
\newcommand{\LL}{\mathcal K} 

\newcommand{\One}{\mathbf{1}}
\newcommand{\dOne}[1]{d\mspace{1mu}\One_{#1}}
\newcommand{\Roots}{\mathrm{Roots}}
\DeclareMathOperator{\im}{im}

\newcommand{\Dmult}[1]{{W_{(s_1, \ldots, s_{#1})}}}

\usepackage{graphicx}   
\newcommand{\WedgeAux}[2]{\scalebox{1.35}{$#1\wedge$}}
\newcommand{\Wedge}{\mathop{\mathpalette\WedgeAux\relax}\nolimits}

\newcommand\KSV{monthly-kutin-smithline-vatter}

\begin{document}

\title[Max-degree filtration and Multicolor Permutation Wordle]
{The max-degree filtration, Koszul complexes, and \\
 dimension counts for Multicolor Permutation Wordle}

\date{\today}

\author[Lawren Smithline]{Lawren~Smithline}
\address{Center for Communications Research, 805 Bunn Drive, Princeton,
NJ 08540-1966, USA}
\email{lawren@ccr-princeton.org}

\hypersetup{pdftitle={Max-degree filtration, Koszul complexes, and Multicolor Permutation Wordle}, pdfauthor={Lawren Smithline}}

\begin{abstract}
    We explore the Koszul resolutions of a family of ideals, each generated by
  multivariate polynomials that split completely into linear factors.
    The motivating examples come from the game \classic of Kutin and Smithline,
  as analyzed by Kutin, Smithline, and Vatter, and an open question,
  Conjecture 4.5 in Kutin and Smithline's
  work on \rainbow.
  For each ideal in this family, with our chosen generators,
  the Koszul differentials are strict for the max-degree filtration and,
  more generally, at every multidegree whose coordinates form a set of
  consecutive integers.
  The resulting dimension count resolves Kutin and Smithline's conjecture.
  
  The proof is an induction in which each inductive step is an iteration. The core argument is human generated. The present proof, developed with AI assistance, separates the 
  iteration of multidegree reductions from the mapping cone induction. This method avoids the need to chase compatibility of lifts up to the top grade of the complex.
\end{abstract}



\maketitle

\section{Introduction}

Kutin and Smithline~\cite{arxiv-kutin-smithline} introduce the game \classic.
  A guesser attempts to recover a setter's hidden permutation of
  the set $[n] = \{1, \ldots, n\}$ in the fewest guesses.  In each round,
  the guesser submits a permutation, and, as in the game Wordle,
  learns which entries are correct.
  The number of permutations that their strategy \basicstrat solves in
  exactly $r$ rounds is the Eulerian number
  $A(n, r-1)$~\cite[Thm.~3.5]{arxiv-kutin-smithline}.
  Kutin, Smithline, and Vatter~\cite{\KSV} allow any word in $[n]^n$
  as a guess, and show that \basicstrat is optimal in a strong sense, which 
  we call \emph{dominance}, following~\cite[\S7]{arxiv-kutin-smithline}: for 
  every $r$, \basicstrat solves as
  many permutations in $S_n$
  within $r$ rounds as any possible strategy.

  Kutin and Smithline~\cite{arxiv-kutin-smithline}
  also generalize \classic to \rainbow.  The game mechanic for \rainbow
  is the same as for \classic.  Each letter $k$ has an additional attribute,
    color, which takes one of $s_k$ values.  Here, we establish a key fact,
  Theorem~\ref{thm:dimcount}, about a certain quotient ring.
  This fact allows us to affirm Conjecture 4.5 of \cite{arxiv-kutin-smithline}.
  \begin{theorem}\label{thm:rainbow}
  \basicstrat is dominant for \rainbow.
  \end{theorem}

\subsection{Dominance for \classic and uniform \rainbow}\label{subsec:classic}

  There are essentially two phases to the dominance proof in
  \cite{\KSV}.
  The first is conversion of the problem to a dimension-counting question,
  an approach described in \cite{babai:linear-algebra-:},
  \cite{guth:polynomial-meth:}, and \cite{matouv-sek:thirty-three-mi:}.
  The overall space, of dimension $n!$, is the quotient of a
  polynomial ring in $n$ variables by a certain ideal.  Restated in
  language amenable to the generalization, the polynomial ring is
    \[R_n = \BC[x_1, \ldots, x_n].\] The ideal is
    \[I_n= \Bigl(\prod_{j = 1}^n (x_j - k): 1 \leq k \leq n\Bigr).\]

    The second phase of the dominance proof is a reduction
    algorithm, the straightening argument given in Allen~\cite[\S2]{allen:the-descent-mon},
    for an element of $R_n$ to a distinguished element of
    its class in $R_n / I_n$.
    This quotient has an explicit basis as a vector space, the
    Garsia--Stanton descent monomials~\cite[\S9]{garsia:group-actions-o:}; see also~\cite{garsia1980combinatorial}.
    Subsets of this basis can be selected by bounding the
    descent number of the permutation associated to a descent monomial.
    These subsets coincide with the subsets selected by
    the \emph{max-degree filtration}, 
    which bounds the maximum degree of any variable in a descent monomial.
    The resulting filtration on the vector space yields the desired
      dimensions, cumulative Eulerian numbers counting the
    permutations in $S_n$ with at most $d$ descents.

In Section~\ref{subsec:unif_rainbow}, we sketch how the dominance of \basicstrat
extends to uniform \rainbow, with all $s_k$ equal to a common value $s$,
without appeal to Koszul resolutions.
For uniform \rainbow, a variation on the straightening argument
of \cite{\KSV} works, 
with Bagno--Biagioli~\cite{bagno2007colored} 
colored-descent monomials in place of descent monomials.
    
\subsection{Outline}
    
    Section~\ref{sec:setup} begins with the
    setup and notation that generalize the
    first phase to \rainbow, including
    a family of ideals $J(M, N) \subseteq R_N$, for $N\geq\nobreak M \geq\nobreak 1$.
    Remark~\ref{rem:translate} specifies
    the connection.
    The rest of Section~\ref{sec:setup} lays out the
    structures for generalizing the
    second phase.
            
    Section~\ref{sec:generators} uses a mapping cone construction 
    to carry the result forward to one more generator.
    Section~\ref{subsec:regular} proves Corollary~\ref{cor:NotZeroDiv}, the one fact
    about the generators of $J(M, N)$ that the induction uses beyond their shape, 
    by an elementary argument based on their construction, 
    \eqref{eq:littleQconstruct} and \eqref{eq:bigQconstruct}.
     Section~\ref{sec:variables} uses an iteration to carry the result 
     forward to one more variable.
      Section~\ref{sec:induction} assembles the previous two sections
     into the main theorem: for each ideal $J(m, n)$
    and each \emph{gap-free} multidegree,
    the filtered component of the Koszul resolution at that multidegree is exact.
    The ideals $J(n, n)$ are applicable to \rainbow.
    
    The induction runs along diagonals of the grid with axes enumerating number of
variables and number of generators.  Section~\ref{subsec:pred} defines
a strictness predicate $S(N, M)$ and an exactness predicate $E(N, M)$.

\begin{figure}[tbp]
\centering
\begin{tikzpicture}[x=1cm,y=1cm,font=\footnotesize,
  st/.style={draw,rounded corners=2pt,align=center,inner sep=3pt},
  imp/.style={->,thick},
  nxt/.style={->,dashed,rounded corners=4pt},
  bas/.style={->,dotted,thick}]
\node[font=\small] at (0,2.6)     {$(N-1,M-1)$};
\node[font=\small] at (6.4,2.6)   {$(N-1,M)$};
\node[font=\small] at (10.75,2.6) {$(N,M)$};
\node[st] (A1) at (0,0)    {$E(N-1,M-1,b,1)$\\ $\wedge\ S(N-1,M-1,b-s_M)$};
\node[st] (B1) at (6.4,0)  {$E(N-1,M,b,1)$};
\node[st] (C1) at (11.5,1.0)  {$S(N,M,a)$};
\node[st] (C2) at (11.5,-1.0) {$E(N,M,a,1)$};
\node[st] (A2) at (0,-3.6)   {$E(N-1,M-1,b,h+1)$\\ $\wedge\ E(N-1,M-1,b-s_M,h)$};
\node[anchor=south] at ([yshift=1pt]A2.north) {grade $h \geq 1$};
\node[st] (B2) at (6.4,-3.6) {$E(N-1,M,b,h+1)$};
\node[st] (C3) at (11.5,-3.6) {$E(N,M,a,h+1)$};
\node[st] (P1) at (6.4,1.75)  {$E(N-1,1,b,1)$};
\node[st] (P2) at (6.4,-1.85) {$E(N-1,1,b,h+1)$};
\node[st] (P0) at (0,-1.75)   {$S(1,1)$};
\draw[imp] (A1) -- node[above]{Lemma~\ref{it:conezero}}
                   node[below]{uses Cor.~\ref{cor:NotZeroDiv}} (B1);
\draw[imp] (A2) -- node[below]{Lemma~\ref{it:coneplus}} (B2);
\draw[bas] (P1) -- node[left,align=right]{base case,\\ Prop.~\ref{prop:Ebase}} (B1);
\draw[bas] (P2) -- node[right,align=left]{base case,\\ Prop.~\ref{prop:Ebase}} (B2);
\draw[bas] (P0) -- node[left,align=right]{base case, \\ Prop.~\ref{prop:Sbase}} (A1);
\draw[imp] (B1.east) -- node[pos=0.9,above left=-1pt,align=center]{Prop.~\ref{prop:reduction}\\ chain $D\,\One_N\to a$} (C1.west);
\draw[imp] (B1.east) -- node[pos=0.85,below left=-5pt,align=center]{Prop.~\ref{prop:Ereduction}\\ chain $a\to-\One_N$} (C2.west);
\draw[imp] (B2) -- node[below,align=center]{Prop.~\ref{prop:Ereduction}\\ chain $a\to-\One_N$} (C3);
\draw[nxt] (C1.north) |- (1.55,3.05) -- ([xshift=15.5mm]A1.north);
\draw[nxt] (C2.east) -- ++(0.35,0) |- (-1.55,3.4) -- ([xshift=-15.5mm]A1.north);
\draw[nxt] (C2.east) -- ++(0.95,0) |- (-1.55,-5.6) -- ([xshift=-15.5mm]A2.south);
\draw[nxt] (C3.east) -- ++(0.2,0) |- (1.55,-5.25) -- ([xshift=15.5mm]A2.south);
\node[fill=white,inner sep=2pt] at (5.5,3.4)  {next step: $N\mapsto N+1$, $M\mapsto M+1$};
\node[fill=white,inner sep=2pt] at (5.5,-5.6) {next step: $N\mapsto N+1$, $M\mapsto M+1$};
\end{tikzpicture}
\caption{The induction of Theorem~\ref{thm:induction}.  
The column headers show (number of variables, number of generators).
Solid arrows are implications within one step, which takes the column $(N-1,M-1)$ to
$(N-1,M)$ by the mapping cone and then to $(N,M)$ by reduction.  
The upper part of the diagram shows the interplay between strictness $S$ in grade~$0$ and
exactness $E$ in grade~$1$.  The lower part shows positive grades; for $h \geq 1$, grades
$h$ and $h+1$ on the left feed into grade $h+1$ in the middle and right.
Dotted arrows are base cases: for
$M=1$, the middle column has $E(N-1, 1)$, 
and at $N=M=2$, in the left column, the strictness base case is $S(1,1)$.  
Dashed arrows show the feedback of the conclusions at $(N,M)$
into the next step, with $N \mapsto N+1, M \mapsto M+1$.
The first draft of this figure was composed by Claude~Opus~5.5 in response to a specific prompt.
It is human-reviewed.
}
\label{fig:inductionmap}
\end{figure}
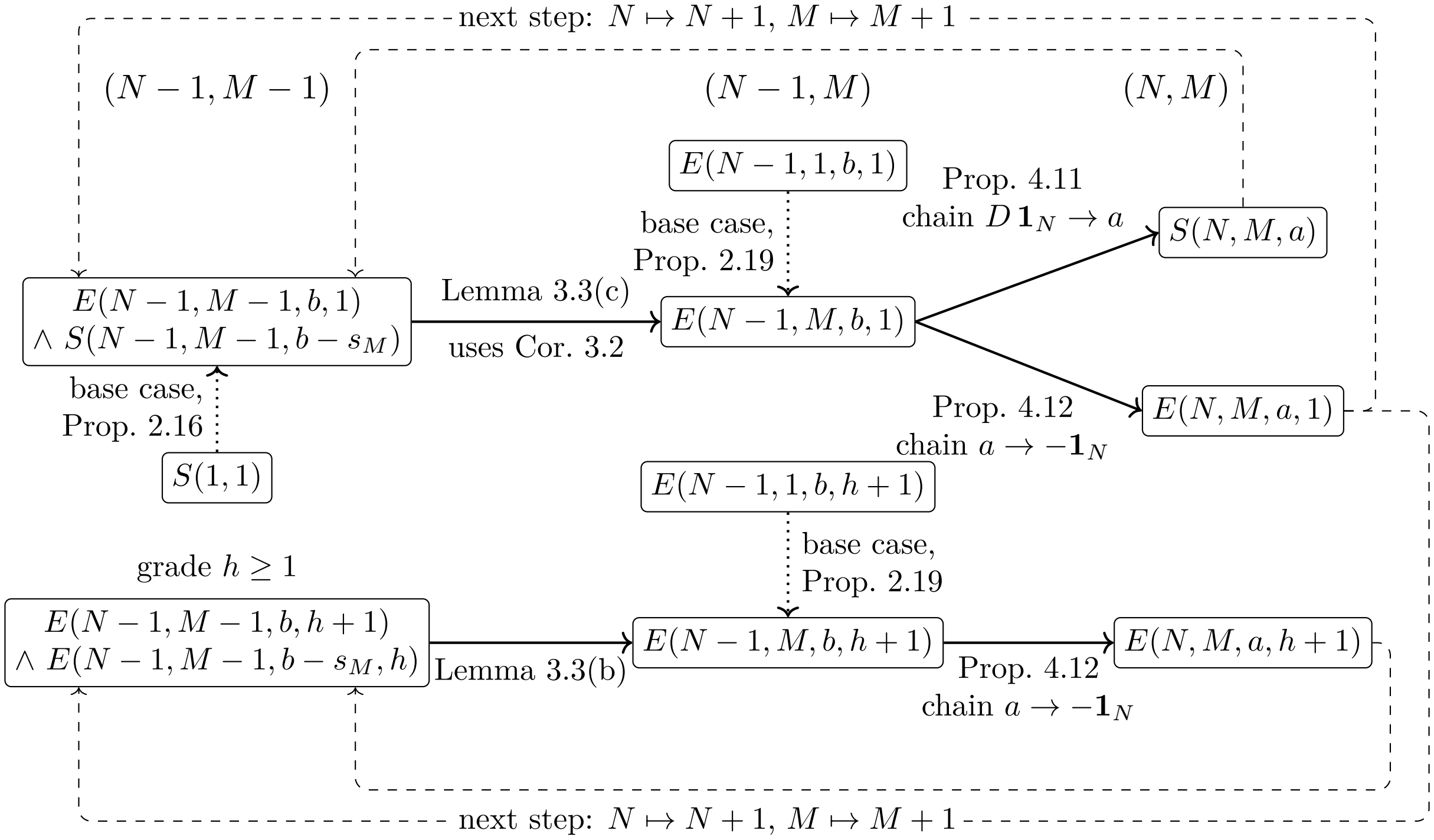

The base cases of Theorem~\ref{thm:induction} are $S(1, 1)$ and $E(N,1)$
(Propositions~\ref{prop:Sbase} and~\ref{prop:Ebase}).
Figure~\ref{fig:inductionmap} traces $S(N, M)$ and $E(N, M)$
through the induction of Theorem~\ref{thm:induction}.

The dimension count result, Theorem~\ref{thm:dimcount} of Section~\ref{sec:dimcount},
uses Theorem~\ref{thm:induction} only through Corollary~\ref{cor:filterexact}.
 Theorem~\ref{thm:dimcount} fills the gap in translating the 
    second phase  of the dominance proof for \classic to \rainbow.
    Section~\ref{subsec:rainbow} explains how to patch Theorem~\ref{thm:dimcount}
    into the argument of \cite{\KSV} to complete the proof of Theorem~\ref{thm:rainbow}.

    Finally, Section~\ref{sec:conclusion} describes the human process and the
    AI collaboration.  The core of the argument, construction
    of a syzygy using an induction wrapped around an iteration, is a human product; 
    both the logical layout and the overall exposition benefit from
    conversations with Anthropic's Claude models Sonnet~5.5, Opus~5.5, and Fable~5.2.

\section{Setup and notation}\label{sec:setup}

For an instance of the Koszul resolution at the entry point of our induction, 
there are $n$ variables and $m \leq n$ generators.
As we construct a syzygy through the course of the proof, we use $N$ and $M \leq N$
as the running parameters of the induction: $N$ variables and $M$ generators.
Where we use $N$ without specifying $n$, it is understood that $N \leq n$.

Fix positive integers $s_1,\ldots, s_n$ and
distinct positive reals $v_1,\ldots, v_n$.
Let
\begin{equation}\label{eq:littleQconstruct}
c_k=v_k^{s_k}, \qquad q_k(x)=x^{s_k}-c_k.
\end{equation}
The key facts we need about these $q_k(x)$ are that each has
simple roots and that any two have disjoint sets of roots.  The binomial shape of these $q_k(x)$
is simply a convenience.

Let $R_N = \BC[x_1, \ldots, x_N]$, and let
\begin{equation}\label{eq:bigQconstruct}
Q(k,N)=\prod_{i=1}^N q_k(x_i) \text{ for } N \leq n.
\end{equation}
\begin{definition1}
For $P \in R_N$, let $\deg_j P$ be the degree of $P$ as a polynomial in $x_j$
with coefficients in the polynomial ring in the $N-1$ variables
$\{x_1, \ldots, x_N\} \setminus \{x_j\}$, with the convention $\deg_j 0 = -\infty$.
\end{definition1}

\begin{definition1}
The \emph{multidegree} of nonzero $P\in R_N$ is the vector
\[ (\deg_1 P, \ldots, \deg_N P ) \in \Z^N. \]
\end{definition1}

Observe that
$\deg_jQ(k,N)=s_k$ for every $j$, and $Q(k,N)=q_k(x_N)\,Q(k,N-1)$.
Since $R_N$ is a domain,
degrees add: $\deg_j(Q(k,N)Z)=s_k+\deg_jZ$.
We abbreviate $Q_k=Q(k,N)$ when $N$ is clear.

Let $J(M, N)$ be the ideal $(Q(k, N) : 1 \leq k \leq M)$, and let $J_N = J(N, N)$.
For \rainbow, the ideal of interest is $J_n$.
For \classic, all $s_k = 1$, all $v_k = k$, and $J_n = I_n$.

\begin{remark1}\label{rem:translate}
For \classic, the first phase of the dominance proof, translation
from permutation guess to recorded polynomial,
involves a factor $x_j - v_k$ when letter $k$ is guessed incorrectly at
position $j$.  For \rainbow, the factor is $x_j - v_k\zeta_{s_k}^{g_k}$,
for $\zeta_{s_k}$ a primitive $s_k$th root of unity,
when letter $k$ with color $g_k$ is guessed incorrectly at position $j$.
\end{remark1}

In the course of the induction, we select one variable at a time for
reduction.  
Each generator $Q_k$ of $J(M, N)$ is invariant under
permutation of the variables $x_1, \ldots, x_N$.
These permutation symmetries induce isomorphic
free resolutions of $J(M, N)$,
described below in Proposition~\ref{prop:permutevar}.

\subsection{Graded components of \texorpdfstring{$\Wedge^* R_N^M$}{the exterior algebra} and the Koszul complex}

Denote the basis of the exterior algebra $\Wedge^* R_N^M$ by $e(T)$ for $T\subseteq[M]$.
The grade of $T$ is $|T|$.

\begin{definition1}
  For $T \subseteq [M]$, 
  the degree $\deg T$ is $\sum_{k\in T} s_k$.
  The degree $\deg e(T)$ of the basis element $e(T)$ is $\deg T$.
\end{definition1}
\begin{definition1}
  For $Z=\sum_TZ(T)e(T)$ define
  \[\nu_j(Z)=\max_{T:\,Z(T)\neq 0}\bigl(\deg_jZ(T)+\deg T\bigr) \]
  with $\nu_j(0) = -\infty$, consistent with $\deg_j 0 = -\infty$.  Define
  \[ \nu (Z) = \bigl( \nu_1(Z), \ldots, \nu_N(Z) \bigr). \]
\end{definition1}
The following notations for multidegrees are helpful: 
\begin{enumerate}
\item for $a, b \in \Z^N$,
write $a \succeq b$ for $a_j \geq b_j$ for $1 \leq j \leq N$;
\item for $a, b \in \Z^N$,
write $a \succ b$ for $a \succeq b$ and $a \neq b$;
\item for $a \in \Z^N$, write $a_{\what p} \in \Z^{N-1}$ for
$a$ with coordinate $p$ deleted and ${\what a}$ for $a_{\what N}$;
\item write $\epsilon_p \in \Z^N$ for the $p$th unit vector of $\Z^N$.
\end{enumerate} 

When we come to Corollary~\ref{cor:filterexact} in
Section~\ref{sec:induction} and dimension counting in
Section~\ref{sec:dimcount}, we shall be interested in multidegrees
that are uniform in all variables, so it is useful to denote by
$\One_N$ the vector with all coordinates equal to 1.

The mapping cone argument of Section~\ref{subsec:cone} 
uses multidegree vectors that differ by a multiple of $\One_N$.
So, for a vector $b \in \Z^N$ and a scalar $s \in \Z$, we write $b-s$ 
as a shorthand for $b - s \One_N$.

Let $\sgn(T,k)=(-1)^{|\{j\in T:\ j<k\}|}$. Note that $\sgn(T\cup\{k\},k)=\sgn(T,k)$.
\begin{definition1} The differential $\delta$ is
\[
\delta\mspace{1mu} e(T)=\sum_{k\in T}\sgn(T,k)\,Q_k\,e(T\setminus \{k\}).
\]
\end{definition1}
\begin{definition1}
  The Koszul complex
  $K(N, M)$ consists of $\Wedge^h R_N^M$
  together with the differential maps
  \[
  \delta:\Wedge^h R_N^M \rightarrow
  \Wedge^{h-1} R_N^M.
  \]
\end{definition1}
We call the index $h$ in the definition of the Koszul
complex the \emph{grade},
rather than homological degree,
to avoid confusion with $\deg_j P$ for $P\in R_N$ and $\deg e(T)$.
Additionally, we sometimes write $\delta_{K(N, M)}$ with the subscript to emphasize
the differential of a particular Koszul complex. 
\begin{remark1}
We state the Koszul differential as a map that lowers the grade by one.
Eisenbud~\cite[\S17.2]{Eisenbud1995} writes the Koszul complex cohomologically,
with a differential that raises the grade; the two versions are isomorphic, and
his $H^j$ corresponds to our grade $M-j$ (cf.\ \cite[Ex.~4.5.2]{Weibel}).
\end{remark1}

\begin{definition1}
 For $b\in \Z^N$, define
\[
F^b_h(N, M) =\{Z\in \Wedge^h R_N^M:\ b \succeq \nu(Z) \}.
\]
\end{definition1}
By degree addition, $\delta(F^b_h(N, M))\subseteq F^b_{h-1}(N, M)$. 
\begin{definition1}
  The filtered Koszul complex $K(N, M, b)$ is the Koszul complex $K(N,M)$
restricted in grade $h$ to $F^b_h(N, M)$.
\end{definition1}
We also denote the grade zero filtered component $F^b_0(N, M)$ by $R_N(b)$,
to emphasize that it is contained in $R_N$ and does not depend on $M$.  If $b_j < 0$
for any $1 \leq j \leq N$, that is, $b \not\succeq 0\, \One_N$,
then $R_N(b) = \{ 0 \}$.

For the reduction in Proposition~\ref{prop:reduction},
it is convenient to start at a multidegree with all coordinates equal.
\begin{definition1}
  The max-degree of $Z$ is $\mu(Z) = \max_{1\leq j \leq N} \nu_j(Z)$.
\end{definition1}
It is a slight abuse of notation to write $F^d_h(N, M)$ for $F^{\dOne{N}}_h(N, M)$
and $R_N(d)$ for $R_N(\dOne{N})$. 

In Sections~\ref{sec:generators}, \ref{sec:variables}, and~\ref{sec:induction},
we use the following notation
to condense the indexing of $K$ and its graded components:
\begin{align*}
\K &= K(N-1, M-1), \quad &  \Kplus& = K(N-1, M), \quad & \K' &= K(N, M), \\
\LL(h) &= \Wedge^h R_{N-1}^{M-1}, \quad & \LL^+(h)& = \Wedge^h R_{N-1}^{M}, \quad &
\LL'(h) &= \Wedge^h R_{N}^{M}.
\end{align*}
We use analogous notation for the filtered complexes and their graded components:
\begin{align*}
\FF(b, *) &= K(N-1, M-1, b), \quad &  \FF^+(b, *)& = K(N-1, M, b), \quad & \FF'(b, *) &= K(N, M, b), \\
\FF(b, h) & = F_h^b(N-1, M-1), \quad & \FF^+(b, h) & = F_h^b(N-1, M), \quad &
\FF'(b, h) & = F_h^b(N, M).
\end{align*}

In Section~\ref{sec:variables}, 
we relate filtered Koszul complexes that differ by
one in the variable count $N$. 
In Lemma~\ref{lem:step-gf}, the reduction step,
we permute the variables as needed to always work on reducing the degree of
the last variable.
\begin{proposition}\label{prop:permutevar}
  Let $\sigma$ be a permutation in the symmetric group $S_N$.
  The Koszul complex $K(N, M)$ is invariant under permutation of the
  variables $x_1, \ldots, x_N$ by $\sigma$.  The filtered Koszul complex
  $K(N, M, b)$ is isomorphic to $K(N, M, \sigma(b))$, where $\sigma$ acts
  on $\Z^N$ by permuting coordinates.
\end{proposition}
\begin{proof} The verification is mechanical.
  With reordered
  variables $x^{\sigma}_k = x_{\sigma(k)}$, we use the multidegree
  vector $\sigma(b)$ with $\sigma(b)_k = b_{\sigma(k)}$.
\end{proof}

\subsection{Gap-free multidegree and predicates
  for the induction and iteration argument}\label{subsec:pred}

We want to show exactness of $K(N, M, \dOne{N})$.  More specifically,
we want the homology $H_h(K(N, M, \dOne{N}))$ to be trivial for $1 \leq h \leq M-1$,
the differential on $F^{\dOne{N}}_M(N, M)$ to be injective,
and the differential on $F^{\dOne{N}}_1(N, M)$ to surject onto
$J(M, N) \cap F^{\dOne{N}}_0(N, M)$.

We use a larger class of multidegrees than scalar multiples of $\One_N$.

\begin{definition1}\label{def:gapfree-gf}
A vector $b\in\Z^N$ is \emph{gap-free} if its set of values
$\{b_1,\dots,b_N\}$ is a set of consecutive integers.  
\end{definition1}

We use \emph{translation closure} to refer to the fact that the shift $b\mapsto b-\One_N$
preserves the gap-free property.  The mapping cone argument
of Section~\ref{subsec:cone}  
uses the consequence that if $b$ is gap-free, then so is $b - s_M$.

\begin{definition1}
For $1 \leq M \leq N$ and $b \in \Z^N$, statement $S(N, M, b)$ is
\[ 
\delta(F^b_1(N, M)) = R_N(b) \cap J(M, N).
\]
Statement $S(N, M)$ is: for all gap-free $b\in \Z^N$, $S(N, M, b)$.
\end{definition1}
\begin{remark1}
By degree addition,
\[ \delta(F^b_{1}(N, M)) \subseteq R_N(b)\cap\delta(\Wedge^{1}R_N^M) \]
always, so $S(N, M, b)$ means that $\delta\colon F^b_{1}(N, M)\to R_N(b)$ is \emph{strict}.
\end{remark1}

\begin{proposition}\label{prop:Sbase}
  For every $N\geq 1$, $S(N, 1)$ holds.
\end{proposition}
\begin{proof}
The module $\Wedge^1 R_N^1$ has a single basis element, which we write as $e_1$.
For
\[
Z=\delta(Z_1 e_1)=Q_1  Z_1\in F^b_0(N, 1),
\]
we have $Z_1 e_1 \in F^b_1(N, 1)$ by degree addition:
\[
\deg_j Z_1=\deg_j Z-s_1\le b_j-s_1.
\]
\end{proof}
\begin{remark1}
  Corollary~\ref{cor:filterexact} and the dimension counting argument in Section~\ref{sec:dimcount}
  use $S(N, M, b)$, the surjectivity of the map $\delta: F^b_1(N, M) \rightarrow J(M, N) \cap R_N(b)$.  Proposition~\ref{prop:reduction} derives $S(N, M, b)$ from $E(N-1, M, *, 1)$ in every case except $S(1, 1)$.  That exception, $S(1, 1)$, is a base case.
\end{remark1}

\begin{definition1}
For $1 \leq M \leq N+1$, $b \in \Z^N$, and grade $h > 0$,
statement $E(N, M, b, h)$ is
\begin{equation}\label{eq:statementE}
X \in \ker \delta: F^b_h(N, M)\rightarrow F^b_{h-1}(N, M)
\Rightarrow X \in \im \delta: F^b_{h+1}(N, M) \rightarrow F^b_h(N, M).
\end{equation}
Statement $E(N, M, *, h)$ is: for all gap-free $b\in \Z^N$, $E(N, M, b, h)$.
Statement $E(N, M)$ is: for $1 \leq h \leq M$, $E(N, M, *, h)$.
\end{definition1}
If $h>M$, then $F^b_h(N, M)=0$, so $E(N, M, b, h)$ holds trivially.
Restating $E(N, M, b, h)$: the homology $H_h\bigl(K(N, M, b)\bigr)$ is trivial,
that is, every cycle in $F^b_{h}(N, M)$ is the boundary of an element of
$F^b_{h+1}(N, M)$.
\begin{proposition}\label{prop:Ebase}
For $N\geq 1$, $b \in \Z^N$, and $h > 0$, $E(N, 1, b, h)$ holds.  In particular, $E(N, 1)$ holds.
\end{proposition}
\begin{proof}
The module $\Wedge^1 R_N^1$ has a single basis element, which we write as $e_1$.
An element of $F^b_1(N, 1)$ is $z e_1$ with $z\in R_N$, and
  \[ \delta(z e_1) =  z\, \delta\mspace{1mu} e_1= z Q_1. \]
  Since $R_N$ is a domain, $z Q_1 = 0$ implies $z = 0$, so the only cycle in grade~$1$ is $0$.
For $h > 1$, $E(N, 1, b, h)$ holds trivially, because $\Wedge^h R^1_N = 0$.
\end{proof}
\begin{remark1}
  The proofs of Propositions~\ref{prop:Sbase} and~\ref{prop:Ebase} do
  not use the gap-free multidegree hypothesis.
\end{remark1}

\section{Incrementing the generator count}\label{sec:generators}

The left side of Figure~\ref{fig:inductionmap}
increments the generator count using a mapping cone construction.
The result of Lemma~\ref{lem:cone} in grade 1 requires an element of
grade~$0$ to be the boundary of an element of a filtered component.
Grade~$0$ has no exactness statement, so this takes two steps.
The non-zero-divisor result of Section~\ref{subsec:regular} shows that
the element is a boundary, and strictness, $S(N-1,M-1)$, shows that it
is the boundary of an element of the filtered component.

\subsection{A regularity result for the ideal generators}\label{subsec:regular}

We show, in Corollary~\ref{cor:NotZeroDiv},
that
$Q^-_M$ is not a zero divisor modulo $J(M-1, N-1)$.
We use a descent that, in its
intermediate steps, reduces the generator set before the variable set,
the opposite order from the descent used for Theorem~\ref{thm:induction}.
Corollary~\ref{cor:NotZeroDiv} is the only place where
Theorem~\ref{thm:induction} relies on the
distinctness of the $v_k$.

A standard proof
uses that a polynomial ring over a field
is Cohen--Macaulay~\cite[Prop.~18.9]{Eisenbud1995}.
We can avoid that theorem, because by
construction, each $Q_k$ factors completely into
linear factors, and for $k \neq \ell$, $Q_k$ and $Q_\ell$ have no
common factor.

To condense the indexing of the generators and ideals, write
\[
Q^-_k = Q(k, N-1), \qquad
\cJ = J(M-1, N-1), \quad \cJ^- = J(M-1, N), \quad \cJ' = J(M, N).
\]
Let $\Roots(q_k)$ be the set of the $s_k$ roots of $q_k(x)$.
By construction~\eqref{eq:littleQconstruct}, the roots of $q_k$ have absolute
value $v_k$, and the $v_k$ are distinct.  So
\begin{equation}\label{eq:qKnonzero}
q_\ell(\beta) \neq 0 \quad \text{whenever } \beta \in \Roots(q_k)
\text{ and } k \neq \ell.
\end{equation}
By~\eqref{eq:bigQconstruct}, $Q_k$ is the product of the linear forms
$x_i - \beta$ with $i \in [N]$ and $\beta \in \Roots(q_k)$.

\begin{proposition}[Linear forms are not zero divisors]\label{prop:linform}
Let $N \geq M \geq 0$, $j \in [N]$, and $\alpha \in \BC$ satisfy
$q_k(\alpha) \neq 0$ for $1 \leq k \leq M$.
If $Y \in  R_N$ and $(x_j - \alpha) Y \in \cJ'$, then $Y \in \cJ'$.
\end{proposition}
\begin{proof}
We argue by induction on $M$ for all $N \geq M$ at once.  
For $M = 0$, $R_{N} / \cJ' = R_{N}$ is a domain.
In this case, $(x_j - \alpha) Y = 0$ implies $Y = 0$.

Suppose the statement is true at $M-1$ and every number of variables $N' \geq M-1$.
Since each $Q_k$ is symmetric in the variables, we may assume $j = N$.
Let $\pi: R_N \rightarrow R_{N-1}$ be the ring homomorphism
that evaluates $x_N$ at $\alpha$.  The kernel of $\pi$ is the ideal $(x_N - \alpha)$.
Since $\pi(x_i) = x_i$ for $1 \leq i < N$, $\pi$ is surjective.
Factoring $Q_k$ as $q_k(x_N)Q^-_k$, we see that
\begin{equation}\label{eq:QkImage}
 \pi\bigl(Q_k\bigr) = q_k(\alpha) Q^-_k, \quad q_k(\alpha) \neq 0 \text{ for } 1 \leq k \leq M.
 \end{equation}
By hypothesis, there are $P'_1, \ldots, P'_M \in R_N$ such that
\begin{equation}\label{eq:idealelement}
(x_N - \alpha) Y = \sum_{k=1}^M P'_k Q_k.
\end{equation}
Write $P_k = q_k(\alpha) \pi(P'_k)$.

\smallskip
\emph{Claim 1:  $P_M \in \cJ$.}   Apply $\pi$ to \eqref{eq:idealelement}.
The left side is in $\ker \pi$.  By \eqref{eq:QkImage},
\[
0 = \sum_{k=1}^M P_k Q^-_k.
\]
Isolating the term $P_M Q^-_M$, we see that $P_M Q^-_M \in \cJ$.
The polynomial $Q^-_M$ is a product of linear forms $x_i - \beta$ with
$i \in [N-1]$ and $\beta \in \Roots(q_M)$.
By~\eqref{eq:qKnonzero},
$q_\ell(\beta) \neq 0$ for $1 \leq \ell < M$.
The induction hypothesis for $M-1$ generators and $N-1$ variables
applies, and we cancel these linear forms one at a time.  Hence, $P_M \in \cJ$.

\smallskip
\emph{Claim 2: $P'_M \in \cJ^- + (x_N - \alpha)$.}
By \eqref{eq:QkImage}, $\pi$ maps the generators $Q_k$ of $\cJ^-$ to
nonzero constant multiples of the generators of $\cJ$, so 
$\pi\bigl(\cJ^-\bigr) = \cJ$.  Since $q_M(\alpha) \neq 0$, Claim 1 shows
$\pi(P'_M) \in \cJ$, so
\[ \exists\, \TP \in \cJ^- :  \pi(\TP) = \pi(P'_M). \]

For this $\TP$, we have
\[ P'_M - \TP \in \ker \pi = (x_N - \alpha). \]
Rewriting,
\begin{equation}\label{eq:YM}
\exists\, Y_M \in R_N : P'_M = \TP + Y_M (x_N-\alpha).
\end{equation}

\smallskip
Using \eqref{eq:YM} to expand the last term of the sum in \eqref{eq:idealelement}, we derive
\[
(x_N - \alpha) \bigl(Y - Y_M Q_M \bigr) = 
\sum_{k=1}^{M-1} P'_k Q_k +  \TP Q_M \in \cJ^-.
\]
By the induction hypothesis for $M-1$ generators and $N$ variables,
$Y - Y_M Q_M \in \cJ^-$.  Hence, $Y \in \cJ'$.
\end{proof}

\begin{corollary}\label{cor:NotZeroDiv}
Let $N \geq M \geq 1$.  If $Y \in R_{N-1}$ and $YQ^-_M \in \cJ$,
then $Y \in \cJ$.  That is, $Q^-_M$ is not a zero divisor modulo $\cJ$.
\end{corollary}
\begin{proof}
The polynomial $Q^-_M$ factors completely into linear forms
$x_j - \alpha$, with $j \in [N-1]$, and $\alpha \in \Roots(q_M)$, and,
by \eqref{eq:qKnonzero},
 $q_k(\alpha) \neq 0$ for $1 \leq k < M$.
Proposition~\ref{prop:linform}, for $M-1$ generators and $N-1$ variables,
cancels these linear forms one at a time,
yielding $Y \in \cJ$.
\end{proof}

\subsection{The mapping cone lemma}\label{subsec:cone}
 
\begin{lemma}\label{lem:cone}
Let $2\le M\le N$, $b\in\Z^{N-1}$, $h > 0$, and $s=s_M$.  If
\begin{thmlist}
\item $E(N-1,M-1,b,h)$;
\item\label{it:coneplus}  if $h \geq 2$, then $E(N-1,M-1,b-s,h-1)$; and
\item\label{it:conezero} if $h = 1$, then $S(N-1,M-1,b-s)$,
\end{thmlist}
then $E(N-1,M,b,h)$.
\end{lemma}

\begin{figure}[tbp]
\[
\begin{tikzcd}[column sep=4.5em, row sep=3.2em]
\FF(b,h+1) \arrow[r,"\delta_\K"] \arrow[d,hook]
& \FF(b,h) \arrow[r,"\delta_\K"] \arrow[d,hook]
  \arrow[from=ddl,"(-1)^hQ^-_M"' {pos=0.3}]
& \FF(b,h-1) \arrow[d,hook]
  \arrow[from=ddl,"(-1)^{h-1}Q^-_M"' {pos=0.3}] \\
\FF^+(b,h+1) \arrow[r,dashed,"\delta_\Kplus" {pos=0.3}, crossing over] \arrow[d,"\pi"]
& \FF^+(b,h) \arrow[r,dashed,"\delta_\Kplus" {pos=0.3}, crossing over] \arrow[d,"\pi"]
& \FF^+(b,h-1) \arrow[d,"\pi"] \\
\FF(b-s,h)\wedge e_M \arrow[r,"{\tilde \delta}_\K"']
& \FF(b-s,h-1)\wedge e_M \arrow[r,"{\tilde \delta}_\K"']
& \FF(b-s,h-2)\wedge e_M
\end{tikzcd}
\]
\caption{The mapping cone of multiplication by $Q^-_M$
restricted to $\FF^+(b, *)$, as in the
proof of Lemma~\ref{lem:cone}.  The columns are grades $h+1$, $h$, and
$h-1$ of the short exact sequence \eqref{eq:exactcone}.  
A piece with a negative grade is zero.
The bottom row has the induced differential $\tilde \delta_\K$ of \eqref{eq:eMdelta}
that makes \eqref{eq:exactcone} an exact sequence of complexes, 
so all of the squares commute.
The arrows in the middle row are dashed because $\delta_\Kplus$ is not
the direct sum of the maps in the top and bottom rows: it is the sum of
the three arrows that leave the top and bottom pieces of a column.
Claim~1 in the proof of Lemma~\ref{lem:cone} is a lift 
along the bottom row, and Claim~2 is a lift 
along the top row.
Claude Opus~5.5 produced the first draft of this figure
in response to a specific prompt.  It is human-reviewed.}
\label{fig:cone}
\end{figure}
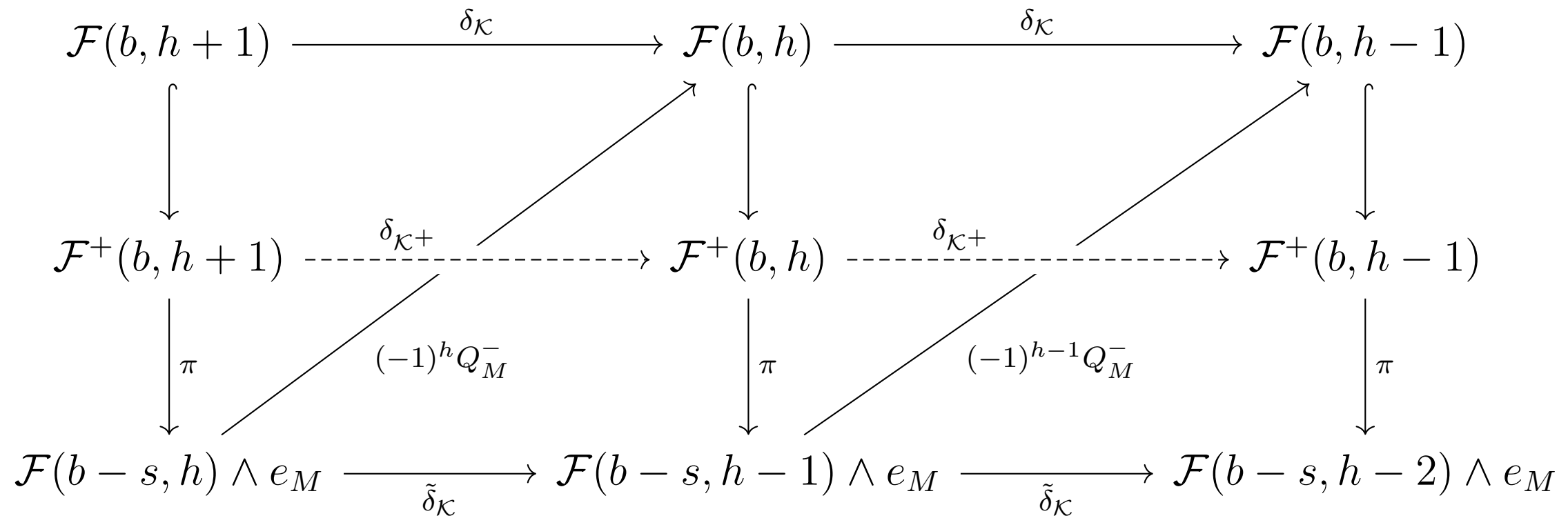

The proof of Lemma~\ref{lem:cone} follows the diagram in Figure~\ref{fig:cone}.
The columns are the exact sequences in grades $h+1$, $h$, and $h-1$ 
for the filtered part $\FF^+(b, *)$ of $\Kplus$ as
the mapping cone of multiplication by $Q^-_M$ on $\K$~\cite[Cor.~17.11]{Eisenbud1995}; see also \cite[\S1.5 and Lemma~4.5.3]{Weibel}:
\begin{equation}\label{eq:exactcone}
0 \to \FF(b, *) \hookrightarrow \FF^+(b, *)\xrightarrow{\ \pi\ }
\FF(b-s, *)\wedge e_M \to 0,
\end{equation}
where $\pi$ is projection to the $e_M$ part of the complex.
The third complex has grade-$j$ component $\FF(b-s,j-1)\wedge e_M$
and differential 
\begin{equation}\label{eq:eMdelta}
{\tilde \delta}_\K(Y\wedge e_M)=(\delta_\K Y)\wedge e_M.
\end{equation}
Starting with a cycle $X$ in the middle row, middle column,
the proof proceeds from \eqref{eq:conesplit} to \eqref{eq:Wtogether}
and produces $W$ such that $\delta W = X$.

\begin{proof}

\smallskip
Let $X$ satisfy the premise of $E(N-1, M, b, h)$:
\[
X \in \FF^+(b, h), \ \delta X = 0.
\]
Rephrasing \eqref{eq:exactcone}, 
  $X \in \LL^+(h)$ has a unique expression
  \begin{equation}\label{eq:conesplit}
  X = X_0 + X_1\wedge e_M \ \text{ for } \
  X_0 \in \LL(h), \
  X_1 \in \LL(h-1). 
  \end{equation}
The filtration operates on each coefficient $X(T)$ separately.  For
$T\subseteq[M-1]$, the coefficient $X(T)$ is $X_0(T)$; and for $T'\subseteq[M-1]$,
the coefficient $X(T'\cup\{M\})$ is $X_1(T')$, with $\deg(T'\cup\{M\})=\deg T'+s$.  So
\begin{equation}\label{eq:eMpart}
 X\in \FF^+(b, h) \iff X_0\in \FF(b,h) \text{ and } X_1\in \FF({b-s},{h-1}).
\end{equation}

\smallskip
 Since $\sgn(T,k)$ indicates the parity of the count of elements of $T$ below $k$, 
 \[
  \sgn(T'\cup\{ M\},k)  =\sgn(T',k) \quad   \text{for } 1 \le k\le M, \text{ and }
  \]
\[
  \sgn(T'\cup\{M\}, M) =\sgn(T',M) =(-1)^{|T'|} \quad \text{for } T' \subseteq [M-1]. \]
For any grade $j$ and $Z \in \LL(j)$,
\begin{equation}\label{eq:conedelta}
\delta(Z\wedge e_{M})= \bigl( (\delta Z)\wedge e_{M}\bigr) +(-1)^j Q^-_M Z.
\end{equation}
Equation~\eqref{eq:conedelta} for $Z = X_1$ and $\delta X=0$ together imply
$\delta X_1 = 0$ for $h \geq 2$ and
\begin{equation}\label{eq:cyc}
\delta X_0+(-1)^{h-1}Q^-_M X_1=0, \quad h \geq 1.
\end{equation}

\smallskip
\emph{Claim 1: There exists $Y_1 \in \FF({b-s}, h)$
such that $X_1 = \delta Y_1$.}
If $h\ge2$, then $X_1$ is a cycle in $\FF(b-s, h-1)$ of positive grade, and
by $E(N-1,M-1,b-s,h-1)$, there exists $Y_1\in\FF(b-s,h)$ such that
$X_1=\delta Y_1$.  If $h=1$, then $X_1\in R_{N-1}$ and \eqref{eq:cyc}
gives $Q^-_MX_1\in\cJ$.  By Corollary~\ref{cor:NotZeroDiv},
$X_1\in\cJ=\delta\,\LL(1)$.  So $X_1$ satisfies the premise of
$S(N-1,M-1,b-s)$, and again, there exists $Y_1\in\FF(b-s,1)$ with
$X_1=\delta Y_1$.

Let \[
W_1 = Y_1 \wedge e_{M} \in \FF^+(b, h+1). \]
By \eqref{eq:conedelta},
\begin{equation}\label{eq:deltaW1}
\delta W_1 = \delta Y_1 \wedge e_{M} + (-1)^{h}Q^-_M Y_1.
\end{equation}

\smallskip
\emph{Claim 2: There exists $W_0 \in \FF(b, h+1)$ such that
$X'= X - \delta W_1 = \delta W_0$.}
Since $\delta X = 0$ and $\delta W_1$ is a boundary, $\delta X' =0$.
Rewriting $X'$ using \eqref{eq:conesplit} and \eqref{eq:deltaW1},
we have
\begin{equation}\label{eq:longXprime}
X' = X_0 + X_1 \wedge e_M -  \delta Y_1 \wedge e_{M} - (-1)^{h}Q^-_M Y_1.
\end{equation}
By Claim 1, $\delta Y_1 = X_1$, so the $\FF({b-s},{h-1})\wedge e_M$ parts of \eqref{eq:longXprime}
cancel.  By degree addition, $Q^-_M Y_1 \in \FF(b, h)$, and
$X_0 \in \FF(b, h)$ by \eqref{eq:eMpart}, so
\begin{equation}\label{eq:shortXprime}
 X' = X_0 -  (-1)^{h}Q^-_M Y_1 \in \FF(b, h)
 \end{equation}
is a cycle. The conclusion of $E(N-1,M-1,b,h)$ gives
$W_0\in\FF(b,h+1)$ with $\delta W_0=X'$.

The constructions of Claim 1 and Claim 2 show
\begin{equation}\label{eq:Wtogether}
W=W_0+W_1\in \FF^+(b, h+1)
\end{equation}
and $\delta W=X$.  We have shown $E(N-1, M, b, h)$.
\end{proof}
\smallskip

\begin{remark1}\label{rem:coneLES}
In the language of homological algebra, \eqref{eq:exactcone}
is a short exact sequence of
complexes with long
exact homology sequence~\cite[\S1.3]{Weibel} containing
\[ 
H_h\bigl(\FF(b,*)\bigr) \to H_h\bigl(\FF^+(b,*)\bigr) \to
H_{h-1}\bigl(\FF(b-s,*)\bigr)
\xrightarrow{\ (-1)^{h-1}Q^-_M\ } H_{h-1}\bigl(\FF(b,*)\bigr).
\] 
By~\eqref{eq:conedelta}, the connecting map is multiplication by $(-1)^{h-1}Q^-_M$.
So $E(N-1,M,b,h)$, the vanishing of 
$H_h\bigl(\FF^+(b,*)\bigr)$, follows from the two claims above.
Claim~1 is the injectivity of the connecting map. In grades
$h\ge2$, the domain of the connecting map, $H_{h-1}\bigl(\FF(b-s, *)\bigr)$, is $0$, by $E(N-1,M-1,b-s,h-1)$.
In grade $h=1$, the domain of the connecting map is 
$R_{N-1}(b-s)/\bigl(\cJ\cap R_{N-1}(b-s)\bigr)$, by $S(N-1,M-1,b-s)$,
and its action, multiplication by $Q^-_M$, is injective by
Corollary~\ref{cor:NotZeroDiv}.  
Claim~2 is $E(N-1,M-1,b,h)$, that is, $H_h\bigl(\FF(b, *)\bigr) = 0$.
\end{remark1}

\begin{proposition}\label{prop:cone}
  For $N \geq M > 1$, if $S(N-1, M-1)$ and $E(N-1,M-1)$, then $E(N-1, M)$.
\end{proposition}
\begin{proof}
Let $b \in \Z^{N-1}$ be gap-free and let $h > 0$.

Gap-free vectors are closed under translation, so $b - s_M$ is gap-free.
The hypotheses, with the trivial cases of grade above $M-1$, include
\begin{thmlist}
\item $E(N-1,M-1,b,h)$; 
\item if $h \geq 2$, then $E(N-1,M-1,b-s_M, h-1)$; and
\item if $h = 1$, then $S(N-1,M-1,b-s_M)$.
\end{thmlist}
By Lemma~\ref{lem:cone}, $E(N-1,M,b,h)$ holds.
\end{proof}

\begin{remark1}\label{rem:cone}
The initial Claude conversation~\cite{claude-sonnet-2026-koszul}
 advised that the decomposition \eqref{eq:conedelta}
  is the standard mapping-cone description of the Koszul complex,
  while
  the bookkeeping of the multidegree filtration
  and the use of $S(N-1, M-1)$ and $E(N-1,M-1)$ in place of ordinary exactness
  are new.

      The mapping-cone description fills a gap in an earlier human-produced draft;
   see Section~\ref{sec:conclusion}.
\end{remark1}

\section{Incrementing the variable count}\label{sec:variables}

The right side of Figure~\ref{fig:inductionmap} increments the variable count.
Section~\ref{subsec:topcoeff} relates filtered Koszul complexes that differ
by one in the variable count, through the top coefficient map and the
exact sequence~\eqref{eq:tcexact}.
Section~\ref{subsec:chain} constructs chains of gap-free multidegrees, and
Section~\ref{subsec:reduction} iterates Lemma~\ref{lem:step-gf} along these chains.

\subsection{The top coefficient map and the Koszul differential}\label{subsec:topcoeff}
As a consequence of Proposition~\ref{prop:permutevar}, when two Koszul
complexes differ by one in the variable count, we may take the extra
variable to be $x_N$.

\begin{definition1}\label{def:tc}
For a multidegree $b\in\Z^N$ and a grade $h$, define the top coefficient map
\[ \tc:F^b_h(N, M) \rightarrow F^{\what b}_h(N-1, M), \qquad
   \tc(Y)=\sum_{|T|=h}\bigl(\text{coefficient of }x_N^{b_N-\deg T}\text{ in }Y(T)\bigr)\,e(T), \]
where a coefficient with $b_N-\deg T<0$ is zero.
\end{definition1}
The coefficient of $x_N^{b_N-\deg T}$ in $Y(T)$ has degree at most
$b_j-\deg T$ in each $x_j$ with $j<N$, so $\tc(Y)\in F^{\what b}_h(N-1, M)$.
\begin{proposition}\label{prop:tcexact}
Let $N\ge2$, $b\in\Z^N$, and $h\ge0$, and suppose $b_N\ge\min_{j<N}b_j$.
The sequence
\begin{equation}\label{eq:tcexact}
0 \to F^{b-\epsilon_N}_h(N, M) \hookrightarrow
 F^b_h(N, M) \xrightarrow{\ \tc\ }   F^{\what b}_h(N-1, M) \to 0
\end{equation}
is a split exact sequence of vector spaces, with section
\begin{equation}\label{eq:tcsection}
\sigma(W)=\sum_{|T|=h}x_N^{b_N-\deg T}\,W(T)\,e(T), 
\quad W \in F^{\what b}_h(N-1, M).
\end{equation}
\end{proposition}
\begin{proof}
For $|T|=h$, the $e(T)$ summands of the three spaces of \eqref{eq:tcexact}
are
\[
R_N(b-\epsilon_N-\deg T)\,e(T),\quad R_N(b-\deg T)\,e(T)\quad\text{and}\quad
R_{N-1}({\what b}-\deg T)\,e(T),
\]
and $\tc$ and $\sigma$ are each
direct sums of linear maps on corresponding summands.  So, it suffices to
analyze $\tc$ and $\sigma$ for a single $T$.

If $b_N-\deg T<0$, then
\[
R_N(b-\epsilon_N-\deg T) = R_N(b-\deg T) = R_{N-1}({\what b}-\deg T) = 0.
\]
The first two are because $b_N - \deg T < 0$ 
and the third is because  $\min_{j<N}b_j \leq b_N < \deg T$. 
Otherwise, as a vector space,
\begin{equation}\label{eq:singlesplit}
R_N(b - \deg T) = \bigoplus_{k = 0}^{b_N - \deg T} R_{N-1}({\what b} - \deg T)\, x_N^k.
\end{equation}
On this summand, $\tc$ is the projection to the top term, $k=b_N-\deg T$;
the kernel of $\tc$ is the sum of the terms with $k<b_N-\deg T$, which is
$R_N(b-\epsilon_N-\deg T)$; and $\sigma$ is the inclusion of the top term,
so the composition ${\tc}\circ\sigma$ is the identity.  This gives the exactness and the
splitting of~\eqref{eq:tcexact} and the section of \eqref{eq:tcsection}.
\end{proof}
\begin{proposition}\label{prop:tcchain}
Let $b\in\Z^N$ and $h \geq 0$.
On $F^b_{h+1}(N, M)$,
\begin{equation}\label{eq:tcDelta}
 \delta_{K(N-1,M)} \tc = \tc \delta_{K(N, M)}.
 \end{equation}
That is, $\tc$ is a map of complexes from $K(N, M, b)$ to $K(N-1, M, {\what b}).$
\end{proposition}
\begin{proof}
Each element of $F^b_{h+1}(N, M)$ is a sum of terms $P\,e(T)$ with $|T|=h+1$
and $P\in R_N(b-\deg T)$.
The maps $\tc$, $\delta_{K(N-1, M)}$, and $\delta_{K(N, M)}$ are linear,
so it suffices to take a single term $P\,e(T)$. 
If $b_N-\deg T<0$, then $P=0$ and both sides of \eqref{eq:tcDelta}
are zero.
Otherwise, let $P^{\tc}\in R_{N-1}$ be the coefficient of
$x_N^{b_N-\deg T}$ in $P$, so that
\[ \tc\bigl(P\,e(T)\bigr)=P^{\tc}e(T). \]
For $k\in T$, $P\,Q(k,N)=P\,q_k(x_N)\,Q(k,N-1)$, and
$\deg(T\setminus\{k\})=\deg T-s_k$.
Since $\deg_N P\le b_N-\deg T$, only the term $x_N^{s_k}$ of
$q_k(x_N) = x_N^{s_k} - c_k$ reaches $x_N^{b_N-\deg(T\setminus\{k\})}$,
and $Q(k, N-1) \in R_{N-1}$, so
\[
\tc \bigl( P\,Q(k,N) e(T\setminus\{k\})\bigr) =  P^{\tc}\,Q(k,N-1) e(T\setminus\{k\}).
\]
Therefore
\begin{align*}
\tc \bigl( \delta_{K(N,M)}\, P\,e(T)\bigr)
 &=\sum_{k\in T}\sgn(T,k)\,P^{\tc}\,Q(k,N-1)\,e(T\setminus\{k\})\\
 &=\delta_{K(N-1,M)}\, P^{\tc}e(T)
  =\delta_{K(N-1,M)}\,\tc\bigl(P\,e(T)\bigr).\qedhere
\end{align*}
\end{proof}
\begin{remark1}
The section $\sigma$ of \eqref{eq:tcsection} is not the unique choice of section
and not a map of complexes.
\end{remark1}

\subsection{Chains of gap-free multidegrees with deletion property}\label{subsec:chain}

The reduction argument of Section~\ref{subsec:reduction}
follows a chain of gap-free multidegrees, which we construct in
Proposition~\ref{prop:chain-gf} below.

\begin{lemma}[Choice of reduction coordinate\footnote{\label{ft:fable}
Claude Fable~5.2 composed the first draft of this lemma, with a proof by contradiction.  We give a direct proof.}]\label{lem:choice-gf}
Let $N\ge2$, and let $a,b\in\Z^N$ be gap-free with $a \succ b$.  
Among the indices $j$ with $a_j>b_j$, let $p$ be
one with $a_p$ largest.  Then
\begin{thmlist}
\item\label{it:choicemin} $a_p\ge\min_{j\ne p}a_j$;
\item\label{it:choicedel} $a_{\what p}$ is gap-free;
\item\label{it:choicestep} $a-\epsilon_p$ is gap-free;
\item\label{it:goodstep} $a-\epsilon_p \succeq b$.
\end{thmlist}
\end{lemma}

\begin{proof}
We show that either $a_p = \max_{j\in[N]} a_j$ or that there is $k \neq p$
such that $a_k = a_p$.

Suppose $(a_p - 1)\One_N \succeq b$.  For $j\neq p$ such that $a_j > b_j$, $a_p \geq a_j$, and for $j\neq p$ such that $b_j \geq a_j$, $a_p > b_j$, so also $a_p \geq a_j$.
Since $N \geq 2$, $[N] \setminus \{p\}$ is not empty, so there is $j \neq p$ such that
$a_j \leq a_p$. This verifies part~(a) and shows that $a_p = \max_{j\in[N]} a_j$.  
Since $a_p$ is maximal among the coordinates of $a$, parts~(b) and~(c) follow: $a_{\what p}$ and $a - \epsilon_p$ are both gap-free.

Alternatively, $(a_p - 1)\One_N \not\succeq b$, so there is some $j \neq p$
such that $b_j \geq a_p$.  Since $b$ is gap-free and $a_p > b_p$, there is $k \neq p$,
possibly $k = j$, such that $b_k =  a_p$.  Since $a \succ b$, $a_k \geq b_k$, and by the
maximality of $a_p$, $b_k \geq a_k$.  So $b_k = a_k = a_p$, and
the coordinate $a_k$ verifies part~(a).  For parts~(b) and~(c), 
$a_k = a_p$ ensures
$a_{\what p}$ is gap-free and $a - \epsilon_p$ is gap-free.

In either case, $a \succ b$ and $a_p > b_p$, so $a_p-1 \geq b_p$ and thus, 
part~(d) holds: $a - \epsilon_p \succeq b$.
\end{proof}

\begin{proposition}\label{prop:chain-gf}
Let $N\ge2$, and let $a,b\in\Z^N$ be gap-free with $a\succeq b$.
Put $r=\sum_{j=1}^N(a_j-b_j)$.  There are gap-free vectors
\[
a=a^{(0)}\succ a^{(1)}\succ\cdots\succ a^{(r)}=b
\]
and indices $p_0,\ldots,p_{r-1}\in[N]$ such that, for $0\le i<r$,
\begin{thmlist}
\item\label{it:chainmin} $a^{(i)}_{p_i}\ge\min_{j\ne p_i}a^{(i)}_j$;
\item\label{it:chaindel} $a^{(i)}_{\what p_i}$ is gap-free;
\item\label{it:chainstep} $a^{(i+1)}=a^{(i)}-\epsilon_{p_i}$.
\end{thmlist}
\end{proposition}

\begin{proof}
We argue by induction on $r$.  If $r=0$, then $a=b$, and the chain has no
steps.  If $r>0$, then $a\succ b$.  Let $p_0$ be the index chosen in
Lemma~\ref{lem:choice-gf}.  Parts \ref{it:choicemin} and~\ref{it:choicedel} give
parts \ref{it:chainmin} and~\ref{it:chaindel} for $i=0$, part~\ref{it:choicestep} shows
that $a^{(1)}=a-\epsilon_{p_0}$ is gap-free, and part~\ref{it:goodstep}
shows $a^{(1)}\succeq b$.
Since $\sum_j(a^{(1)}_j-b_j)=r-1$, the induction hypothesis gives the rest
of the chain.
\end{proof}

\begin{remark1}\label{rem:chain-gf}
The chain of Proposition~\ref{prop:chain-gf} is the sequence of variable
choices in the proof of Proposition~\ref{prop:reduction}: at step $i$, the
variable $x_{p_i}$ is reduced by one degree, by Lemma~\ref{lem:step-gf}.
Conditions~\ref{it:chainmin}--\ref{it:chainstep} supply what
that lemma needs: \ref{it:chainmin} is its hypothesis~\ref{it:min-gf};
\ref{it:chaindel} makes its hypothesis~\ref{it:lift-gf} available from
$E(N-1, M)$; and \ref{it:chainstep}, with $a^{(i+1)}\succeq\nobreak b$, gives its
hypothesis~\ref{it:room-gf}.
\end{remark1}

\subsection{The reduction step}\label{subsec:reduction}

Proposition~\ref{prop:reduction}, $E(N-1, M, *, 1) \Rightarrow S(N, M)$ for $N\geq 2$ and $N \geq M \geq 1$,
and  Proposition~\ref{prop:Ereduction}, $E(N-1, M) \Rightarrow E(N, M)$ for $N \geq M \geq 2$, 
both use an iteration of Lemma~\ref{lem:step-gf}.

\begin{lemma}[Reduction by one degree for one variable]\label{lem:step-gf}
Let $N\ge2$, $N \geq M \geq 1$, and $0\le h<M$.  Let $a,b\in\Z^N$ and $p\in[N]$.  
Assume
\begin{thmlist}
\item\label{it:min-gf} $a_p\ge\min_{j\ne p}a_j$;
\item\label{it:lift-gf} $E(N-1,M,a_{\what p},h+1)$ holds;
\item\label{it:room-gf} $a_p>b_p$.
\end{thmlist}
If $\TZ\in \FF'(a, h+1)$ and $\delta\TZ\in \FF'(b, h)$, then there is $\widehat W\in \FF'(a, h+2)$ such that
\[ \TZ-\delta\widehat W \in \FF'(a-\epsilon_p, h+1). \]
\end{lemma}

\begin{figure}[tbp]
\[
\begin{tikzcd}[column sep=large, row sep=large]
\FF'(a-\epsilon_N,h+2) \arrow[r,"\delta"] \arrow[d,hook]
& \FF'(a-\epsilon_N,h+1) \arrow[r,"\delta"] \arrow[d,hook]
& \FF'(a-\epsilon_N,h) \arrow[d,hook] \\
\FF'(a,h+2) \arrow[r,"\delta"] \arrow[d,"\tc"']
& \FF'(a,h+1) \arrow[r,"\delta"] \arrow[d,"\tc"]
& \FF'(a,h) \arrow[d,"\tc"] \\
\FF^+(\what a,h+2) \arrow[r,"\delta"]
  \arrow[u, dashed, start anchor=west, end anchor=west, bend left=70,
         "\sigma: W\mapsto\widehat W"]
& \FF^+(\what a,h+1) \arrow[r,"\delta"]
& \FF^+(\what a,h)
\end{tikzcd}
\]
\caption{The top coefficient map $\tc$ in grades $h+2$, $h+1$, and $h$, in
the proof of Lemma~\ref{lem:step-gf}.  The dashed arrow is not a map of
complexes.
Claude Opus~5.5 produced this figure in response to a specific prompt.
It is human-reviewed.}
\label{fig:tc}
\end{figure}
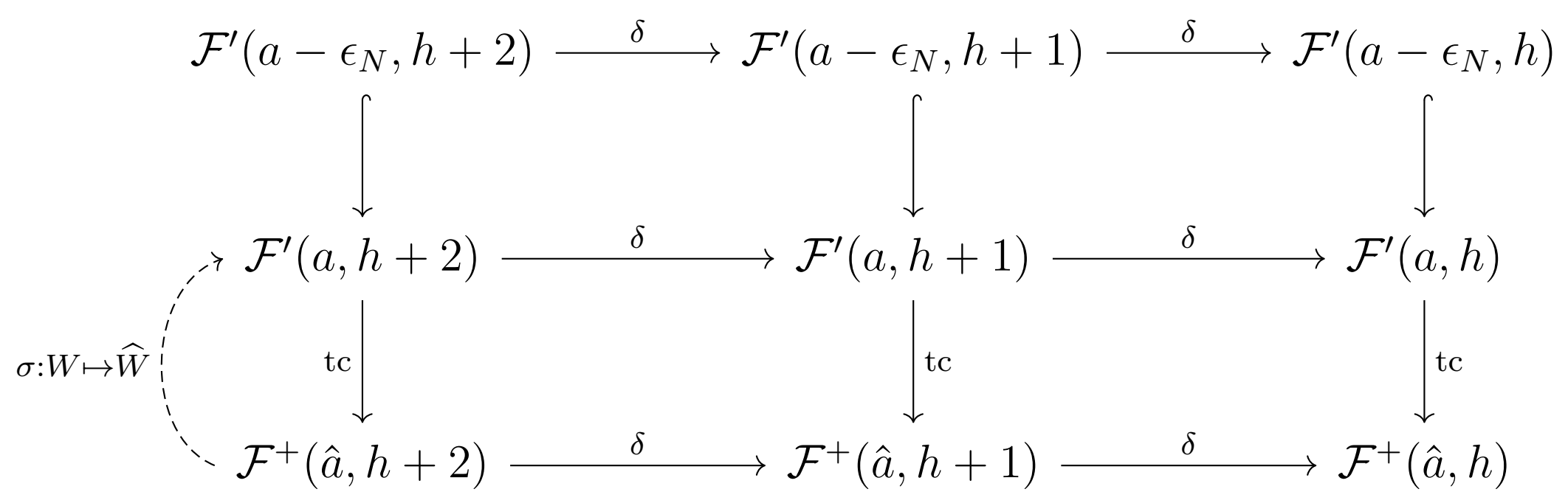

The proof of Lemma~\ref{lem:step-gf}
follows the diagram in Figure~\ref{fig:tc}.
With $p = N$ as in the proof, 
the rows are the complexes $\FF'(a-\epsilon_N,*)$,
$\FF'(a,*)$, and $\FF^+(\what a,*)$, and the columns are the grades $h+2$,
$h+1$, and $h$.
The proof starts with $\TZ$ in the middle row, middle column,
and proceeds step by step, \eqref{eq:Zlift} to \eqref{eq:tzPrime}.

\begin{proof} 
By Proposition~\ref{prop:permutevar}, we may cyclically shift variables $x_p,\ldots,x_N$,
moving $x_p$ to the last place and each of $x_{p+1},\ldots,x_N$ one place down,
shift the coordinates of vectors $a$ and $b$ in the same way,
 and write the reduction in terms of $x_N$.  
Let  $Z=\delta\TZ$.

\smallskip
\noindent\emph{Showing $\tc(\TZ)$ is a cycle.}
By degree addition, $Z=\delta\TZ\in\FF'(a,h)$.
By~\ref{it:room-gf}, $a_N > b_N$, so 
\begin{equation}\label{eq:Zlift}
Z\in\FF'(a,h)\cap\FF'(b,h) \subseteq \FF'(a-\epsilon_N,h).
\end{equation}
By~\ref{it:min-gf} and Proposition~\ref{prop:tcexact} in grade $h$,
\begin{equation}\label{eq:Zkernel}
\tc(Z) = 0.
\end{equation}
By Definition~\ref{def:tc}, $\tc(\TZ)\in \FF^+({\what a}, h+1)$.
By Proposition~\ref{prop:tcchain},
\begin{equation}\label{eq:tcTZzero}
 \delta\,\tc(\TZ)=\tc(\delta\TZ)=0.
\end{equation}

\smallskip
\noindent\emph{The lift and correction.}
By~\ref{it:lift-gf}, $E(N-1,M,{\what a},h+1)$,
\begin{equation}\label{eq:liftexist}
\exists\, W\in \FF^+({\what a}, h+2) : \delta W=\tc(\TZ).
\end{equation}
By~\ref{it:min-gf}, Proposition~\ref{prop:tcexact} applies to $a$ in grade
$h+2$.  With $\sigma$ as in~\eqref{eq:tcsection}, let 
\begin{equation}\label{eq:Wsection}
\widehat W=\sigma(W)\in\FF'(a,h+2),
\end{equation}
so that $\tc(\widehat W)=W$.
By Proposition~\ref{prop:tcchain},
\begin{equation}\label{eq:lowerleft}
\tc(\delta_{\K'}\widehat W)=\delta\,\tc(\widehat W)=\delta W=\tc(\TZ), 
\end{equation}
so
\begin{equation}\label{eq:ZdelWkernel}
\tc(\TZ - \delta_{\K'} \widehat W) = 0.
\end{equation}
By~\ref{it:min-gf}
and Proposition~\ref{prop:tcexact} in grade $h+1$,
\begin{equation}\label{eq:tzPrime}
\TZ' = \TZ-\delta_{\K'} \widehat W\in \FF'(a-\epsilon_N, h+1).
\end{equation}

Another application of Proposition~\ref{prop:permutevar} restores the variables
and vectors $a$ and $b$ to their original order, 
making the reduction effective for $x_p$.  
\end{proof}

\begin{remark1}
Condition~\ref{it:min-gf}, which holds whenever $a_p$ is the largest coordinate of
$a$, makes $\tc$ surjective, by Proposition~\ref{prop:tcexact}.
Within the proof of Lemma~\ref{lem:step-gf}, 
the existence of the lift $W$ using condition~\ref{it:lift-gf}, at \eqref{eq:liftexist},
is the only place where a class of multidegrees enters. 
\end{remark1}

\begin{remark1}\label{rem:stepLES}
We can restate Lemma~\ref{lem:step-gf} in the language of homological algebra,
using $p = N$ for convenience.
By condition~\ref{it:min-gf} and Proposition~\ref{prop:tcexact},
the columns of Figure~\ref{fig:tc} are part of
a short exact sequence of complexes
\[ 0\to \FF'(a-\epsilon_N,*) \xhookrightarrow{\ \iota\ } \FF'(a,*)
   \xrightarrow{\ \tc\ } \FF^+(\what a,*)\to 0,\]
where $\iota$ is the inclusion.
Condition~\ref{it:lift-gf} says $H_{h+1}\bigl(\FF^+(\what a,*)\bigr)=0$, so
the long exact homology sequence~\cite[\S1.3]{Weibel} contains
\[
H_{h+1}\bigl(\FF'(a-\epsilon_N,*)\bigr) \xrightarrow{\ \iota_*\ } H_{h+1}\bigl(\FF'(a,*)\bigr)
   \to 0 \to H_h\bigl(\FF'(a-\epsilon_N,*)\bigr) \xrightarrow{\ \iota_*\ } H_h\bigl(\FF'(a,*)\bigr). \]
The map $\iota_*$ in grade $h+1$ is surjective: a cycle in $\FF'(a,h+1)$ differs by a
boundary from a cycle in $\FF'(a-\epsilon_N,h+1)$.  Along a chain of multidegrees, this is
Proposition~\ref{prop:Ereduction}.  The map $\iota_*$ in grade $h$ is injective: a cycle in
$\FF'(a-\epsilon_N,h)$ that is a boundary from $\FF'(a,h+1)$ is a boundary
from $\FF'(a-\epsilon_N,h+1)$.  For $h=0$, along a chain of multidegrees, this is
Proposition~\ref{prop:reduction}.  Lemma~\ref{lem:step-gf} proves one
statement that yields both: an element of $\FF'(a,h+1)$ whose boundary lies in
$\FF'(a-\epsilon_N,h)$ differs by a boundary from an element of
$\FF'(a-\epsilon_N,h+1)$.  Put differently, the relative homology
\[
H_{h+1}\bigl(\FF'(a, *), \FF'(a - \epsilon_N, *) \bigr) \cong H_{h+1}\bigl(\FF^+({\what a}, *) \bigr)
\]
is zero, which is condition~\ref{it:lift-gf}.
\end{remark1}

\begin{proposition}\label{prop:reduction}
For $N \geq 2$ and $N \geq M \geq 1$, if $E(N-1, M, *, 1)$, then $S(N, M)$.
\end{proposition}

\begin{proof}
Let $b \in \Z^N$ be gap-free 
and let $Z=\delta\TZ\in R_N(b)$ with $\TZ\in\LL'(1)$.
We assume $Z \ne 0$, otherwise there is nothing to prove.
Let $D=\max\bigl( \mu(\TZ), \max_j b_j \bigr)$ and
$a=D\, \One_{N}$.  
The vector $a$ is gap-free, $a \succeq b$, and
$\TZ$ lies in $\FF'(a, 1)$.
Take the chain $a=a^{(0)}\succ\cdots\succ a^{(r)}=b$ and the indices
$p_0,\ldots,p_{r-1}$ of Proposition~\ref{prop:chain-gf}, and put $\TZ^{(0)}=\TZ$.

\smallskip
Let $0\le i<r$, and suppose that
$\TZ^{(i)}\in \FF'(a^{(i)}, 1)$ with $\delta\TZ^{(i)}=Z$.
The hypotheses of Lemma~\ref{lem:step-gf} hold
for $a^{(i)}$, $b$, $p=p_i$, and $h=0$: hypothesis~\ref{it:min-gf} by
Proposition~\ref{it:chainmin}; hypothesis~\ref{it:lift-gf} by
Proposition~\ref{it:chaindel} and the hypothesis $E(N-1,M,*,1)$ of the present
proposition; and hypothesis~\ref{it:room-gf}
because $a^{(i)}_{p_i}=a^{(i+1)}_{p_i}+1>b_{p_i}$, by
Proposition~\ref{it:chainstep}.
Lemma~\ref{lem:step-gf} gives $\widehat W_i\in\FF'(a^{(i)}, 2)$ such that
$\TZ^{(i+1)}=\TZ^{(i)}-\delta\widehat W_i$ lies in $\FF'(a^{(i+1)}, 1)$.
The iterated boundary $\delta\delta \widehat W_i$ is zero, so
$\delta\TZ^{(i+1)}=\delta\TZ^{(i)}=Z$.

After $r$ steps, $\TZ^{(r)}\in \FF'(b, 1)$ and $\delta\TZ^{(r)}=Z$.
\end{proof}

\begin{proposition}\label{prop:Ereduction}
For $N \geq M \geq 2$, if $E(N-1, M)$, then $E(N, M)$.
\end{proposition}
\begin{proof}
Let $a\in\Z^N$ be gap-free, let $0\le h < M$, and let $\TZ\in\FF'(a,h+1)$
with $\delta\TZ=0$.  We show $\TZ\in\delta\,\FF'(a,h+2)$.

If $a$ has a negative coordinate, then $\FF'(a,h+1)=0$ and $\TZ=0$.
Otherwise, let $b=-\One_N$;  $b$ is gap-free, $a\succeq b$, and
$\FF'(b,h+1)=0$, because $b$ has a negative coordinate.
Proposition~\ref{prop:chain-gf} produces the chain $a=a^{(0)}\succ\cdots\succ a^{(r)}=b$
and the indices $p_0,\ldots,p_{r-1}$.
Let $\TZ^{(0)}=\TZ$.

Let $0\le i<r$, and suppose that $\TZ^{(i)}\in\FF'(a^{(i)},h+1)$ is a cycle
with $\TZ-\TZ^{(i)}\in\delta\,\FF'(a,h+2)$.  Its boundary, zero, lies in
$\FF'(b,h)$, so the hypotheses of Lemma~\ref{lem:step-gf} hold for
$a^{(i)}$, $b$, $p=p_i$, and the grade $h$, as in the proof of
Proposition~\ref{prop:reduction}, with $E(N-1,M,*,h+1)$ in place of $E(N-1,M,*,1)$.  The lemma gives
$\widehat W_i\in\FF'(a^{(i)},h+2)\subseteq\FF'(a,h+2)$ such that
$\TZ^{(i+1)}=\TZ^{(i)}-\delta\widehat W_i$ lies in $\FF'(a^{(i+1)},h+1)$.  
Since $\TZ^{(i)}$ and $\delta\widehat W_i$ are both cycles, $\TZ^{(i+1)}$
is a cycle, and $\TZ-\TZ^{(i+1)}\in\delta\,\FF'(a,h+2)$.

After $r$ steps, $\TZ^{(r)}\in\FF'(b,h+1)=0$, so $\TZ\in\delta\,\FF'(a,h+2)$.
\end{proof}

\begin{remark1}\label{rem:original}
The initial Claude conversation~\cite{claude-sonnet-2026-koszul} stated, in
paraphrase, that 
it knew of no reference for Section~\ref{subsec:reduction} 
and believed that the iterated reduction was the new part of the argument.
While Lemma~\ref{lem:step-gf} echoes the reduction step of
Buchberger's algorithm~\cite[Ch.~2, \S6, Thm.~6, and 
\S10, Thm.~6]{cox-little-oshea2025},
the gap-free multidegree, nested iteration,
and mapping cone construction are different from Gr\"obner basis reduction.
\end{remark1}

\begin{remark1}\label{rem:fanout} 
The lift $E(N-1, M)$ used in each iteration is itself obtained from
$S(N-1, M-1)$ and $E(N-1,M-1)$ by Proposition~\ref{prop:cone}, and these are
proved by the same iterations one variable and one generator down.
Unwinding the induction gives a nested iteration, a ``fan-out'' over
variables, in which each layer re-applies the symmetry of the
variables.  Separating the predicates $S$ and $E$ keeps every layer
self-contained: Proposition~\ref{prop:cone} adds a generator once, for all
subsets $T$ simultaneously, and the iterations of
Propositions~\ref{prop:reduction} and~\ref{prop:Ereduction} add a variable.

Because Lemma~\ref{lem:cone} works on all subsets $T$ at once, nothing has to be
made compatible.  Reducing each coordinate $X(T)$
separately requires choosing a single $\tilde z\in\LL^+(h+2)$
shared by all $T\subsetneq T'$ with $|T|=|T'|-1$.  Independently chosen local corrections need not
come from one $\tilde z$.
\end{remark1}

\section{The induction argument and the exactness of filtered complexes}\label{sec:induction}

\begin{theorem}\label{thm:induction}
  For all $N \geq M \ge 1$, $S(N, M)$ and $E(N,M)$ hold: for every gap-free
  $b\in\Z^N$,
  \begin{thmlist}
  \item an element of $\delta\, \K'(1)$ lying in $R_N(b)$ is
  the boundary of an element of $\FF'(b, 1)$;
  \item for every $h \geq 1$, a cycle of $\FF'(b,h)$ is
  the boundary of an element of $\FF'(b,h+1)$.
  \end{thmlist}
  In particular, this holds for every instance of the $n$-variable problem
  with $m \le n$ generators $Q_1,\dots,Q_m$.
\end{theorem}
\begin{proof}
Map the parameter $N$ to $\Delta = N-M \geq 0$.
For particular $\Delta$, we proceed by induction on $M$.
For $M=1$, $E(1+\Delta,1)$ is Proposition~\ref{prop:Ebase}.
For $\Delta=0$, $S(1,1)$ is Proposition~\ref{prop:Sbase}; for $\Delta\geq1$,
Proposition~\ref{prop:reduction} with $M=1$ gives $S(1+\Delta,1)$ from
$E(\Delta,1,*,1)$.

Let $M\ge2$.
We have $S(\Delta+M-1, M-1)$ and $E(\Delta+M-1,M-1)$ by the induction
hypothesis.
Proposition~\ref{prop:cone} gives $E(\Delta+M-1,M)$.  
With the addition of $E(\Delta+M-1,M)$ to the induction hypothesis,
Proposition~\ref{prop:reduction} gives $S(\Delta+M,M)$, and
Proposition~\ref{prop:Ereduction} gives $E(\Delta+M,M)$.
\end{proof}

\begin{remark1}\label{rem:Mzero}
The induction can also start at $M=0$, with base cases $S(N,0)$ and $E(N,0)$,
which hold trivially.  That start needs conventions for $N=0$ and $M=0$:
$R_0=\BC$, each $Q(k,0)=1$ as an empty product,
$F^b_h(0,M)=\Wedge^h R_0^M$ for the empty vector $b$,
and $\Wedge^0 R_N^0=R_N$, so that $F^b_0(N,0)=R_N(b)$.  It also needs
condition~\ref{it:min-gf} to be vacuous for $N=1$, and separate treatment
of the case $N=1$ in either Lemma~\ref{lem:step-gf} or
Proposition~\ref{prop:tcexact}.
We start at $M=1$ to avoid these exceptions.
\end{remark1}

\begin{corollary}[Exactness of the filtered complexes]\label{cor:filterexact}
  Let $N \geq M \ge1$ and $b\in \Z^N$ be gap-free.
The filtered subcomplex of $\K'$,
\[
0\to \FF'(b, M) \to\cdots\to \FF'(b, 1)\to \FF'(b, 0) =R_N(b),
\]
is exact, except at the right end, where
\[
\delta\, \FF'(b, 1) = J(M, N) \cap R_N(b)=\sum_{k\le M}Q_k\,R_N(b-s_k).
\]
\end{corollary}
\begin{proof}
For $1\le h\le M$, exactness at $\FF'(b,h)$ is $E(N,M,b,h)$, part of
Theorem~\ref{thm:induction}.
At the right end, $S(N,M, b)$ gives
  \[
  J(M, N) \cap R_N(b)=\bigl( \delta\LL'(1) \bigr)\cap \FF'(b, 0)
  \subseteq\delta\, \FF'(b, 1),
  \]
  and the
  reverse inclusion is clear.
\end{proof}

\begin{remark1}\label{rem:Szero}
Strictness is used only in grade~$0$.  Theorem~\ref{thm:induction},
through Lemma~\ref{lem:cone}, uses $S(N-1,M-1,b-s)$, and
Corollary~\ref{cor:filterexact} uses $S(N,M,b)$ to describe
$J(M,N)\cap R_N(b)$.  There is no exactness statement in grade~$0$, since
every element of grade~$0$ is a cycle.  Apart from those grade-$0$
statements, the route to $E(N,M)$ passes only through exactness
statements.  Strictness in a grade $h > 0$ follows from
$E(N,M,b,h)$, since a boundary lying in $\FF'(b,h)$ is a cycle lying in
$\FF'(b,h)$.
\end{remark1}

\begin{remark1}
Claude Fable~5.2 was prompted for Remarks~\ref{rem:gapneeded} and~\ref{rem:uniform}.

Remark~\ref{rem:gapneeded} states that the gap-free
hypothesis of Theorem~\ref{thm:induction} cannot be dropped,
and it cannot be relaxed to allow gaps of size up to $\min_k s_k$.

Remark~\ref{rem:uniform} states that when all $s_k$ equal a common value $s$, 
$K(N, M, b)$ is a direct sum over $s^N$ filtered Koszul complexes for the single-color case.

Claude Fable~5.2 observed, in Corollary~\ref{cor:koszulexact},
that Theorem~\ref{thm:induction}, with Proposition~\ref{prop:cone} for
$M=N+1$, implies exactness of the Koszul complex $\K'$ except at the right end.
For $M=N+1$ and $N=n$, that corollary needs a generator $Q_{n+1}$, that is,
data $s_{n+1}$ and $v_{n+1}$ beyond those fixed in Section~\ref{sec:setup}.

In an earlier draft, we started with that exactness, cited to \cite[Cor.~17.5]{Eisenbud1995},
as a consequence of regularity.
\end{remark1}
\begin{remark1}
We can go further than Remark~\ref{rem:uniform}.  Let $t = \gcd(s_1, \ldots, s_M)$.
Each $Q_k$ is a polynomial in the powers $x_i^{t}$, so, in the same way,
the filtered Koszul complex $K(N, M, b)$ is a direct sum over $t^N$ filtered
Koszul complexes for the exponents $s_k / t$, at the multidegrees
$\lfloor (b-r)/t\rfloor$ for $r\in\{0,\ldots,t-1\}^N$.  So $S(N,M,b)$ and
$E(N,M,b,h)$ hold whenever each of these multidegrees is gap-free or has a
negative coordinate, a larger class than the gap-free multidegrees.
\end{remark1}

\section{Dimensions of the quotients \texorpdfstring{$R_n(d) / (J_n\cap R_n(d))$}{R\_n(d)/(J\_n cap R\_n(d))}}\label{sec:dimcount}

We apply Corollary~\ref{cor:filterexact} to the question of dominance of \basicstrat for \rainbow.
For fixed $n$, data $s_1, \ldots, s_n$, and $d \in \Z_{\geq 0}$,
following \cite{arxiv-kutin-smithline},
we let $\Dmult{n}(n, d)$ be the number of $(s_1, \ldots, s_n)$-colored permutations
of length $n$ solved by \basicstrat in exactly $d+1$ rounds.
By \cite[Theorem~4.7]{arxiv-kutin-smithline}, the generating function for
$\Dmult{n}(n, d)$ satisfies
\begin{equation}\label{eq:multi-gf}
\sum_d \Dmult{n}(n,d) u^d = (1-u)\prod_k(1-u^{s_k}) \sum_d (d+1)^n u^d. 
\end{equation}
When all $s_k = 1$, $\Dmult{n}(n,d)$ is the Eulerian number $A(n,d)$ and \eqref{eq:multi-gf}
specializes to the Eulerian polynomial $A_n(u)$ and its well-known generating function formulation
\begin{equation}\label{eq:eulerian-gf}
  A_n(u) = \sum_{d=0}^{n-1}A(n, d) u^d =
  (1-u)^{n+1}\sum_{d=0}^{\infty} (d+1)^n u^d.
\end{equation}

\begin{theorem}\label{thm:dimcount}
  For a positive integer $n$,
  positive integers $s_1, \ldots, s_n$, distinct positive reals $v_1, \ldots, v_n$,
  let $J = J_n$.  For every $d \in \Z_{\geq 0}$, 
  the dimensions $\dim R_n(d) / (J\cap R_n(d))$
  satisfy
  \[
  \dim\tfrac{R_n(d)}{J\cap R_n(d)}-\dim\tfrac{R_n(d-1)}{J\cap
    R_n(d-1)} =\Dmult{n}(n,d).
    \]
\end{theorem}

\begin{proof}
  Let $d \geq 0$.
  Corollary~\ref{cor:filterexact} with $N=M=n$, $b = \dOne{n}$, shows
  $K(n,n,\dOne{n})$ is a finite exact complex of finite-dimensional spaces resolving
$R_n(d)/(J\cap R_n(d))$. Since
$\dim F^b_h(n,n)=\sum_{|T|=h} \max(d-\deg T+1, 0)^{\,n}$, we have
\[
\dim R_n(d)/(J\cap R_n(d))=
\sum_{T\subseteq[n]}(-1)^{|T|}\,\max(d-\deg T+1, 0)^{\,n}.
\]
We find
\[
(1-u)^{n+1}\sum_{d\ge0}\dim\frac{R_n(d)}{J\cap R_n(d)}\, u^d = A_n(u) \prod_{k=1}^n(1-u^{s_k}).
\]
Hence, the successive differences have generating function
\begin{equation}\label{eq:dimgen}
\sum_{d\ge0}\Bigl(\dim\tfrac{R_n(d)}{J\cap
  R_n(d)}-\dim\tfrac{R_n(d-1)}{J\cap
  R_n(d-1)}\Bigr)u^d=A_n(u)\prod_{k=1}^n [s_k]_u, \qquad [s]_u = \sum_{j=0}^{s-1} u^j.
\end{equation}

Substituting the rightmost expression of \eqref{eq:eulerian-gf}
for $A_n(u)$ in \eqref{eq:dimgen}, we produce the right side of
\eqref{eq:multi-gf}.  We read off the left side of
\eqref{eq:multi-gf} to get the desired result.
\end{proof}

The right side of \eqref{eq:dimgen}
has degree $(n-1)+\sum(s_k-1)=\sum s_k-1$ and value $n!\prod s_k$ at
$u=1$. So the dimension reaches $n!\prod s_k$ at $d=\sum s_k-1$ and
stays there.

Claude Sonnet~5.5~\cite{claude-sonnet-2026-koszul} selected and worked out the
following examples of dimension counts, which we have verified:
$n=3$, $s=(1,1,1)$: $1,5,6,6,\dots$; and $n=2$, $s=(3,2)$: $1,4,8,11,12,12,\dots$.
Claude Sonnet~5.5 wrote further that these examples matched direct
rank computations for the evaluation of $R_n(d)$ on the zero set $V(J_n)$,
a set of $n!\prod s_k$ points.

\subsection{Review of \basicstrat dominance for \classic with application to \rainbow}\label{subsec:rainbow}

Section~\ref{subsec:classic} reviews the translation of \basicstrat to dimension
counting via the polynomial method sufficiently to motivate
Theorem~\ref{thm:dimcount}.  Here, we review the results of 
\cite{\KSV} needed to prove Theorem~\ref{thm:rainbow},
\basicstrat dominance for \rainbow.
While we note specific differences for \rainbow as compared with \classic, we
direct the reader to \cite{\KSV} for complete proofs.

\begin{proof}[Proof highlights for Theorem~\ref{thm:rainbow}]
We may assume the setter chooses a secret uniformly at random~\cite[\S2]{\KSV}.  As in~\cite[Theorem~11]{\KSV},
we may restrict attention to deterministic strategies.

Given $n$, $s_1, \ldots, s_n$,
for \rainbow, the setter chooses a permutation and coloring $(\sigma, \kappa)$.
At round $t$, the guesser guesses $(\gamma_t, \chi_t)$.  The response is a set of 
correct positions
\[U_t = \{ i: \sigma(i) = \gamma_t(i) \text{ and } \chi_t(\gamma_t(i)) = \kappa(\gamma_t(i))\}.
\]
For a game instance, if $U_d = [n]$, we define all subsequent guesses $(\gamma_t, \chi_t)$ to be
equal to $(\gamma_d, \chi_d)$.

Given a strategy, a sequence $U_1, \ldots, U_d$ is a \emph{feasible history}
if there exists $(\sigma, \kappa)$ that produces it.  For each feasible history, fix a choice
of witness $(\sigma_U, \kappa_U)$ to its feasibility.
The number of feasible histories of length $d$ bounds the number of 
colored permutations solved within $d+1$ rounds: each feasible history yields
at most one solved colored permutation at the next round.

Select $v_1, \ldots, v_n$ distinct positive real values, yielding the ideal $J_n$ and
its zero set $V(J_n)$.   For $1 \leq k \leq n$, write the 
$s_k$ elements of $\Roots(q_k)$
as $r(k, \ell)$ for $0 \leq \ell < s_k$.
For a feasible history $U = U_1, \ldots, U_d$, resulting from
guesses $\gamma_1, \ldots, \gamma_d$,
the history polynomials are
\[ P_{U, t} = \prod_{i \notin U_t} 
\Bigl( x_i - r\bigl(\gamma_t(i), \chi_t( \gamma_t (i) ) \bigr) \Bigr), \qquad
P_U = \prod_{t = 1}^d P_{U, t}.
\]

For any $(\sigma, \kappa)$ consistent with the history $U$, the factor $x_i - r\bigl(\gamma_t(i), \chi_t( \gamma_t (i))\bigr)$ is nonzero
at $x_i = r\bigl(\sigma(i), \kappa( \sigma(i))\bigr)$ because for $i \notin U_t$, 
guess $(\gamma_t, \chi_t)$
differs from $(\sigma, \kappa)$ at position $i$.
Thus, the triangularity
argument of \cite[Proposition~14]{\KSV} shows
the polynomials $\{P_U : U \text{ a feasible history of length } d \}$
are linearly independent as functions on $V(J_n)$.
Furthermore, the history polynomial $P_U$ of a history of length $d$ lies in the
filtered component $R_n(d)$, because it is the product of $d$ polynomials in $R_n(1)$.

So, $\dim \langle P_U : U \text{ a feasible history of length } d \rangle$ is bounded by
$\dim R_n(d) / \bigl(J_n \cap R_n(d) \bigr)$, and so is the number of colored
permutations that any deterministic strategy solves within $d+1$ rounds.
By Theorem~\ref{thm:dimcount}
and \cite[Def.~4.6 and Theorem~4.7]{arxiv-kutin-smithline}, this dimension is
$\sum_{k\le d} \Dmult{n}(n,k)$, the number of colored permutations that \basicstrat
solves within $d+1$ rounds,
concluding the proof of Theorem~\ref{thm:rainbow}.
\end{proof}

\subsection{Sketch of dominance proof for uniform \rainbow}\label{subsec:unif_rainbow}

        The dominance proof of \cite{\KSV} generalizes to uniform \rainbow,
    with all $s_k$ equal to a common value $s$, without appeal to Koszul
    resolutions.  In brief, the colored permutations of $[n]$ with $s$ colors
    have a group structure, the wreath product $G(s, n) = \Z/s\Z \wr S_n$.

    Use the monomial order by sorted degree vector, then lexicographic order,
    as in \cite{\KSV}.
    The reduction of \cite{\KSV} finds the Garsia--Stanton
    descent monomial that divides the top monomial,
    with the variables ordered so that the descent monomial and its cofactor
    both have weakly decreasing exponent sequences.
    Instead, find the Bagno--Biagioli~\cite{bagno2007colored} colored-descent
    monomial and the perfect $s$-power cofactor whose product is the top monomial.
    The induced symmetric polynomial of the $s$-power cofactor is $G(s, n)$-invariant,
    so it differs by a constant from an element of $J_n$, the analog of $I_n$.
     
The top exponent of the colored-descent monomial of $g\in G(s,n)$ is
    $\mathrm{fdes}(g)$, and the generating function of
    \cite[Theorem~A.1]{bagno2007colored} for their $\mathrm{fdes}$, at their $q = 1$
    with their $t$ as our $u$, agrees with \eqref{eq:multi-gf}, which matches
    the performance of \basicstrat shown in
    \cite[Theorem~4.7]{arxiv-kutin-smithline}.
    
    Since Theorem~\ref{thm:rainbow} proves the more general result, 
    we omit the remaining details of this construction for uniform \rainbow.
    
\section{Remarks on AI collaboration}\label{sec:conclusion}

At the time of this writing, in late 2026, most of us are new to
AI collaboration for technical writing, with experience ranging from perhaps two years
for users closer to the AI development world to less than a year for
individuals with modest subscription budgets.
AI workflows and AI risk management continue to evolve.

The results here, building on \cite{\KSV} and \cite{arxiv-kutin-smithline},
began with the intuition that the ideal $J_n$ would be a good generalization of
$I_n$.  At least the first phase of the translation of the polynomial method
would work there.  The importance of $Q_1, \ldots, Q_n$ being
a regular sequence was also understood: in the later writing, the AI insisted that
it could not be waved off as a standard result.

There were a number of false starts on Corollary~\ref{cor:filterexact} and Theorem~\ref{thm:dimcount}.  Mostly, these involved getting stuck with
an obstruction in one variable.  The idea of building a syzygy together
with analysis of the $n=2$ case inspired the iterated reduction.
Progressing to $n=3$ showed that something more was needed.

The first sketch on paper for the general $n$ noted a reduction in variable degree.
Remark~\ref{rem:fanout} mentions the strategy: reduce each coordinate $X(T)$
separately.  
As this proposed induction proceeded up the grades from $e(T)$, it would
effectively build $K(n, m, \dOne{n})$ for $m = |T|$.

That sketch suggested there might be choices for the next variable to reduce
after descending to one fewer generator.
It further asserted that 
different choices would encounter an obstruction 
in the syzygy module at the next higher grade,
recursively up to the top grade, where the differential is injective.
The details were never worked out; the success of the present argument
suggests that the earlier sketch works.  
Sections~\ref{sec:generators} and~\ref{sec:variables} overtook
the need for completing the argument in the sketch.

At this point, it was decided that this particular proof was a good candidate to
study with AI collaboration.  There was good intuition that the method of proof would
succeed.  The balance leaned against working out the
details by hand and in favor of verifying that
those details had been worked out by another.

The paper sketch was transcribed to a simple \LaTeX\ document.
Claude Sonnet 5.5~\cite{claude-sonnet-2026-koszul} responded with its
\LaTeX\ document.  The Claude Sonnet 5.5 argument differs from the human
argument.  Claude Sonnet 5.5 proposed
the separation of the predicates $E$ and $S$, whose role
Remark~\ref{rem:fanout} explains, as well as use of
a family of multidegrees larger than $\Z\, \One_N$.
(The Sonnet~5.5-proposed family was $\{ b\in \Z^N : \max_j b_j - \min_j b_j \leq 1 \}$;
the present writing uses the gap-free multidegrees proposed by Fable~5.2.)

The draft \cite{claude-sonnet-2026-koszul} was the genesis of the present work,
human-evaluated by writing it longhand, line by line.  There was a subsequent round of
mathematical editing in collaboration with Claude Sonnet 5.5,
including restriction of Theorem~\ref{thm:induction} to the main diagonal $M = N$.
Continuing collaboration with Claude Opus 5.5
included less sweeping mathematical editing,
review of human-added context in Section~\ref{subsec:rainbow},
clarifying the connection between Theorem~\ref{thm:rainbow} here
and the optimality proof of \cite{\KSV}; and standard editing tasks, using
prompts like ``workshopping open item ...'' and a proposed replacement paragraph.
 Claude Fable 5.2 commented that restoring
 Theorem~\ref{thm:induction} to cover $M \leq N$ would eliminate the entire section
 devoted to extending the main-diagonal argument to $M \leq N$. 
 The gap-free class enabled a further
 simplification of the reduction iteration, bringing the proof even closer
 to the original handwritten sketch.  The single degree reduction now exists as
 Lemma~\ref{lem:step-gf}, rather than being entwined in the proof of Proposition~\ref{prop:reduction}.
 Asked whether exactness of the layers of the filtration could replace the
 theorem on regular sequences, Claude Fable~5.2 proposed carrying the exactness
 predicate $E(N, M)$ through the induction beside strictness predicates in all grades.
 With that change, Corollary~\ref{cor:koszulexact} proves the exactness of the Koszul complexes of
 the generators of the ideals $J(M, N)$, rather than quoting the standard result
 for regular sequences.  Additionally, the strictness predicate was dropped except
 in grade 0, where $S(N-1, M-1)$ is required to prove $E(N-1, M, *, 1)$ and Corollary~\ref{cor:filterexact}.
 After this point,
 there was a series of closely coached modifications, adding some of the
 language of homological algebra.  In particular, adding Propositions~\ref{prop:tcexact}  and~\ref{prop:tcchain},
 inspired by the homological algebra proof composed by Claude Fable~5.2,
 transformed the proof of Lemma~\ref{lem:step-gf}
 into a chase through Figure~\ref{fig:tc}.

The principal technique for ensuring human responsibility for the resulting draft
was an overall instruction to the AI that it was to write a review document proposing
revisions to the source (the \texttt{.tex} file)
and that it was not allowed to edit the source itself.
The human user maintained the source separately and applied edits, sometimes
as recommended by the AI, and sometimes not.  

On request, the AI wrote companion source files.
For example, Claude Sonnet 5.5 was asked to review a previous version of
Section~\ref{subsec:regular}, the regularity of 
 the sequence $Q_1, \ldots, Q_N$, and produce a proof that did not rely on the Cohen--Macaulay
 property of $R_N$.  Claude Fable 5.2 was asked to shorten the Claude Sonnet 5.5 work.
 The present Section~\ref{subsec:regular} is a human rewrite of the Claude Fable 5.2 plan.

The separation of roles---human editing guided by AI review, and AI review guided
  by human prompts---may be more conservative than necessary.  It does provide the human author a good basis for claiming one hundred percent human review of the text of this paper.

At a high level, the AI models used, Claude Sonnet~5.5, Claude Opus~5.5,
and Claude Fable~5.2, are, together, the wisest reference librarian one ever knew, 
with access to the largest library; the most diligent and
caffeinated graduate student; the sharpest technical editor; and sometimes
the most inattentive high school sophomore---and they are none of these.
Rather, they are multimodal large language models, used here with tools for reading
files, running code, and searching the web.

\section*{Acknowledgments}
    This work would be incomplete without acknowledging the contributions of
    Samuel Kutin to the \classic and \rainbow work~\cite{arxiv-kutin-smithline} and 
    of both Samuel Kutin and Vince Vatter to the \classic work~\cite{\KSV}.
    Kutin first noticed the connection between \basicstrat and 
    Eulerian numbers, which was critical for proving the performance of
    \basicstrat in all cases.
    Vatter brought the polynomial method to our attention and proved the
        optimality of \basicstrat for the original \classic.

\appendix

\section{Claude Fable~5.2 composed text}\label{app:fable-include}

We include the following Claude Fable~5.2 text, as human-reviewed and easily verified.

\begin{remark1}[on Theorem~\ref{thm:induction} and Corollary~\ref{cor:filterexact}]\label{rem:gapneeded}
The gap-free hypothesis of Theorem~\ref{thm:induction} cannot be dropped,
and it cannot be relaxed to allow gaps of size up to $\min_k s_k$.
Let $N = M = 2$.  If $s_1 = s_2 = 1$, then
\[
(x_2 - c_2)\,Q_1 - (x_2 - c_1)\,Q_2 = (c_2 - c_1)(x_2 - c_1)(x_2 - c_2)
\]
is a nonzero element of $J(2, 2) \cap R_2(b)$ for $b = (0, 2)$, while
$R_2(b - 1) = \{0\}$.  So $S(2, 2, (0, 2))$ fails.
If $s_1 = s_2 = 2$, the same identity with $x_i^2$ in place of $x_i$,
multiplied by $x_1$, shows that
\[
Z = x_1 (x_2^2 - c_1)(x_2^2 - c_2)
\]
lies in $J(2, 2) \cap R_2(b)$ for $b = (2, 4)$.  The elements of
$R_2(b - 2) = R_2\bigl((0, 2)\bigr)$ are polynomials in $x_2$ alone, so every
element of $Q_1 R_2(b-2) + Q_2 R_2(b-2)$ is even in $x_1$, and $Z$ is odd
in $x_1$.  So $S(2, 2, (2, 4))$ fails, although the gap of $b$ is
$\min_k s_k$.
\end{remark1}

\begin{remark1}\label{rem:uniform}
Suppose all $s_k$ equal a common value $s$.  The case $s = 1$
determines the multidegrees $b$ for which $S(N, M, b)$ holds.
Each $Q_k$ is a polynomial in $y_i = x_i^s$, $1 \leq i \leq N$, so
\[
R_N = \bigoplus_{r} x^r\, \BC[y_1, \ldots, y_N], \qquad
r \in \{0, \ldots, s-1\}^N, \quad x^r = x_1^{r_1}\cdots x_N^{r_N},
\]
and $\delta$ preserves the corresponding decomposition of $\Wedge^h R_N^M$.
A term $x^r y^\gamma e(T)$ lies in $F^b_h(N, M)$ exactly when
$r_j + s\gamma_j + s|T| \leq b_j$ for all $j$, that is, when
\[
\gamma_j + |T| \leq \beta(r)_j, \qquad
\beta(r)_j = \left\lfloor \frac{b_j - r_j}{s} \right\rfloor .
\]
Thus $K(N, M, b)$ is the direct sum, over $r$, of copies of the filtered
complexes with multidegree $\beta(r)$ for the generators
$\prod_i (y_i - c_k)$, and
$S(N, M, b)$ holds if and only if the corresponding statement for
these generators holds at $\beta(r)$ for every $r$.
By Theorem~\ref{thm:induction}, this is so whenever every $\beta(r)$ is
gap-free or has a negative coordinate.  For example, if
$b_j \equiv s - 1 \pmod{s}$ for all $j$, then every $\beta(r)$ equals
$\bigl(b - (s-1)\bigr)/s$, and $b$ may have gaps of size $s$.
In Remark~\ref{rem:gapneeded}, the multidegree $b = (2, 4)$ with $s = 2$
has $\beta\bigl((1, 0)\bigr) = (0, 2)$.
\end{remark1}

\begin{corollary}[of Theorem~\ref{thm:induction}]\label{cor:koszulexact}
For $N\ge1$ and $1\le M\le N+1$, the Koszul complex $K(N,M)$ of
$Q_1,\dots,Q_M$ over $R_N$ has trivial homology in every positive grade.
\end{corollary}
\begin{proof}
Let $X$ be a cycle of $K(N,M)$ in grade $h\ge1$, and put $d=\mu(X)$, so
that $X\in F^{\dOne{N}}_h(N,M)$.  The vector $\dOne{N}$ is gap-free.  If
$M\le N$, then $E(N,M,\dOne{N},h)$ holds by Theorem~\ref{thm:induction}.
If $M=N+1$, it holds by Proposition~\ref{prop:cone} with $N+1$ in place of
both $N$ and $M$, applied to $S(N,N)$ and $E(N,N)$ from
Theorem~\ref{thm:induction}.  So $X$ is a boundary.
\end{proof}

\end{document}